\documentclass[10pt,a4paper]{article}
\pdfoutput=1 

\usepackage{lineno,hyperref}
\modulolinenumbers[5]

\usepackage{amssymb}
\usepackage{amsmath}
\usepackage{algorithm}     
\usepackage{algpseudocode} 
\usepackage{subfig}
\usepackage{csquotes}
\usepackage{tikz}
\usetikzlibrary{arrows}
\usepackage{mathtools}
\usepackage{bbm}
\usepackage{amsthm}
\usepackage{pgfplots}
\usepackage{float}
\usepgfplotslibrary{external}
\DeclareRobustCommand{\tikzcaption}[1]{\tikzset{external/export next=false}#1}
\DeclareRobustCommand{\tikzref}[1]{\tikzcaption{\resizebox{!}{\refsize}{\ref{#1}}}}

\newcommand{\email}[1]{\hspace*{\stretch{1}}\emph{\texttt{#1}}}

\makeatletter
\def\blfootnote{\xdef\@thefnmark{$\star$}\@footnotetext}
\makeatother
\newenvironment{Authors}%
  {\begin{center}\begin{bfseries}}%
  {\end{bfseries}\end{center}}
\newenvironment{Addresses}%
  {\begin{flushleft}\begin{itshape}}%
  {\end{itshape}\end{flushleft}}

 \usepackage{fancyhdr}  
  
\fancypagestyle{plain}{
	\fancyhead{}
	\fancyhead[C]{\hfill submitted to Elsevier, October 2020}
}
  
\newtheorem{theorem1}{Theorem}[section]
\newtheorem{theorem2}{Theorem}[section]

\newtheorem{lemma}[theorem1]{Lemma}
\newtheorem{remark}[theorem2]{Remark}

  \newcommand{\vertiii}[1]{{\left\vert\kern-0.25ex\left\vert\kern-0.25ex\left\vert #1 
    \right\vert\kern-0.25ex\right\vert\kern-0.25ex\right\vert}}

\begin{document}

\thispagestyle{plain}

\title{A discretize-then-map approach for the treatment of parameterized geometries in model order reduction}
 \date{}
 
 \maketitle
\vspace{-50pt} 
 
\begin{Authors}
Tommaso Taddei$^{1}$,
Lei Zhang$^{1}$
\end{Authors}

\begin{Addresses}
$^1$
IMB, UMR 5251, Univ. Bordeaux;  33400, Talence, France.
Inria Bordeaux Sud-Ouest, Team MEMPHIS;  33400, Talence, France, \email{tommaso.taddei@inria.fr,lei.a.zhang@inria.fr} \\
\end{Addresses}

\begin{abstract}
We present a   general approach for the treatment of parameterized geometries in projection-based model order reduction. 
During the offline stage,
given
(i)  a family of parameterized domains  $\{  \Omega_{\mu}: \mu \in \mathcal{P} \} \subset \mathbb{R}^D$ where $\mu  \in \mathcal{P} \subset \mathbb{R}^P$ denotes a vector of  parameters, 
(ii) a parameterized mapping ${\Phi}_{\mu}$ between a reference domain $\Omega$ and the parameter-dependent domain $\Omega_{\mu}$, and
(iii) a finite element triangulation of  $\Omega$,
we resort to an empirical quadrature procedure to select a subset of the elements of the grid.  During the online stage, we first use the mapping to ``move" the nodes of the selected elements and then we use standard element-wise residual evaluation routines to evaluate  the residual and possibly its Jacobian.
We discuss how to devise an online-efficient reduced-order model and we  discuss the differences with the more standard ``map-then-discretize" approach (e.g., Rozza, Huynh, Patera, ACME, 2007); in particular, we show how the discretize-then-map framework  greatly simplifies the implementation of the reduced-order model.
We apply our approach to a two-dimensional potential flow problem past a parameterized airfoil, and to  the two-dimensional RANS simulations of the flow past the Ahmed body.
\end{abstract}

 

\section{Introduction}
\label{sec:introduction}

\subsection{Treatment of parameterized geometries in model order reduction}

In { s}cience and { e}ngineering, we are often interested in performing parametric studies to assess the influence of geometry on the solution to a given partial differential equation (PDE) of interest.
We  denote by $\{  \Omega_{\mu}: \mu \in \mathcal{P} \}$
a family of parameterized domains  in $\mathbb{R}^D$, $D=2,3$;
then, we introduce a reference domain $\Omega$  and  the parameterized mapping ${\Phi}: \Omega \times \mathcal{P} \to \mathbb{R}^D$ such that ${\Phi}_{\mu}:= {\Phi}(\cdot; \mu)$ is a bijection in $\Omega$ and ${\Phi}_{\mu}(\Omega) = \Omega_{\mu}$.
In this paper, we present a general approach for the treatment of parameterized geometries in projection-based 
parameterized model order reduction (pMOR). The approach exploits the knowledge of the  mapping  ${\Phi}$: 
automatic construction of rapidly-computable mappings for parameterized geometries is beyond the scope of the present work;
{
we refer to the model reduction literature (e.g., 
\cite{lassila2014model,rozza2021basic}) for a thorough discussion.}

To make the setting mathematically precise, we shall introduce a model problem considered in the numerical results.
Given the family of domains $\{  \Omega_{\mu}: \mu \in \mathcal{P} \}$,
we define the Dirichlet datum
$h: \mathcal{P} \to H^{1/2}(\partial \Omega_{\mu})$; then, we seek  
$u: \mathcal{P} \to H^{1}( \Omega_{\mu})$ such that
$u_{\mu}\big|_{ \partial \Omega_{\mu}  } = h_{\mu}$ and 
\begin{equation}
\label{eq:model_problem}
\int_{\Omega_{\mu}} \nabla    u_{\mu}  \cdot \nabla v \; d{x} = 0 
\quad
\forall \, v \in H_0^1( \Omega_{\mu}  ).
\end{equation}
Since the domain $\Omega_{\mu}$ depends on the parameter, pMOR techniques cannot be directly applied to devise a 
reduced-order model (ROM) for \eqref{eq:model_problem}
{ --- 
pMOR techniques are indeed designed to approximate fields over a parameter-independent domain}.

To address this issue, several authors (e.g., \cite{lassila2014model,rozza2021basic,ballarin2016fast,ballarin2017numerical,rozza2007reduced}) have followed a \emph{map-then-discretize} (MtD) approach. First, 
given the parameterized mapping ${\Phi}$, 
we introduce the weak formulation of the mapped problem for $\widetilde{u}_{\mu} = u_{\mu} \circ {\Phi}_{\mu} \in H^1(\Omega)$,
\begin{equation}
\label{eq:mapped_problem}
\int_{\Omega} {K}_{\mu} \nabla \tilde{u}_{\mu} \cdot \nabla v \, d {x} \, = \,  0 \quad
\forall \, v \in H_0^1(\Omega),
\qquad
\tilde{u}_{\mu}|_{\partial \Omega} =\tilde{h}_{\mu},
\end{equation}
where 
${K}_{\mu} = {\rm det} ( \nabla  {\Phi}_{\mu} ) \, 
\nabla  {\Phi}_{\mu}^{-1}  
\nabla  {\Phi}_{\mu}^{-T}$ and
$\tilde{h}_{\mu} =h_{\mu} \circ  {\Phi}_{\mu}$.
Second, we discretize \eqref{eq:mapped_problem} using a standard high-fidelity (HF) solver --- e.g., the finite element (FE) method --- and we 
apply traditional pMOR techniques to construct a ROM. 

Although the MtD approach has been successfully applied to a broad class of PDEs, it suffers from two major limitations.
First, since we effectively modify the variational formulation  of the problem, the MtD approach might require the implementation of new HF routines:
for the problem at hand, the HF solver for the mapped problem should  deal with non-constant coefficients and also  with  terms of the form
$\int_{\Omega} \partial_{x_i} u  \; \partial_{x_j} v \, d {x}$ with $i\neq j$ that are not present in the original problem. Although this does not represent a major limitation for this specific test case, the need for implementing  additional functionalities in the HF code might prevent the application of projection-based pMOR techniques to more challenging problems: {
in section \ref{sec:simplicity_implementation}, we show through a simple example why the MtD approach might be considerably more involved to apply for stabilized FE formulations.}
Second, it might be convenient to choose parameterized  mappings that are  piecewise smooth in $\Omega$
(see the numerical examples in \cite{rozza2007reduced} and also the two model problems of this paper):
in this case, the   underlying mesh used to discretize \eqref{eq:mapped_problem} should be conforming with the coarse-grained partition associated with the definition of the mapping  to ensure convergence of the HF solver at optimal { polynomial} rate.
{ Discontinuities of ${K}_{\mu}$ inside elements of the mesh pose indeed major issues for the HF quadrature rule and ultimately affect  performance of the solver (cf. section \ref{sec:methods_DtM}).} 

\subsection{Objective of the paper and outline}

In this paper,  we present a \emph{discretize-then-map} (DtM) approach to deal with parameterized geometries for projection-based pMOR. 
During the offline stage, we resort to an empirical quadrature (EQ) procedure to optimally choose a low-dimensional subset of the grid elements;
during the online stage, we first use the mapping to ``move''  the nodes of the selected elements and then we resort to  standard element-wise residual evaluation routines to  assemble the reduced system.
We show that our approach  does not require  the HF mesh to be conforming with the coarse-grained partition; furthermore, it does not modify the equations and thus does not require  modifications to  the HF code.

Our point of departure is a high-order isoparametric continuous    FE discretization: we shall consider both Galerkin and Petrov-Galerkin FE formulations.
As explained in section \ref{sec:methods_DtM}, given a mesh over the domain $\Omega$, FE fields (and thus reduced-order spaces) should be interpreted as pairs of a FE vector and a mapping ${\Phi}$ (cf. Equation \eqref{eq:vector2field}): the vector contains the ``coefficients'' of the solution with respect to the FE basis; the mapping defines the deformation of the nodes of the reference mesh.
We define inner products and norms in the reference configuration, while we define parameter-dependent variational forms in the parameter-dependent (physical) configuration. To make our discussion mathematically precise, we   explicitly work with the algebraic representation of the FE fields: depending on the choice of ${\Phi}$, the latter can be linked to the solution to \eqref{eq:model_problem}  and to the solution to \eqref{eq:mapped_problem}.

We discuss in detail implementation aspects that are key to integrate the ROM with available HF codes; in particular, we show that the DtM framework greatly simplifies  the implementation of ROMs, particularly for nonlinear PDEs. For completeness, we also discuss how to enforce Dirichlet boundary conditions (BCs): more in detail, we discuss the use of lifting functions --- also known as ``control function method'' ---  to deal with inhomogeneous Dirichlet BCs.
We refer to the recent work \cite{star2019extension} and to the references therein for a thorough survey on the imposition of Dirichlet BCs in (Petrov-)Galerkin ROMs:  we envision that our
DtM  approach can cope with both strong and weak enforcement of Dirichlet BCs.

The paper is organized as follows.
{ In section \ref{sec:methods_DtM}, we introduce 
the main features of the DtM framework; furthermore, we compare the approach with the more standard MtD framework and we offer a number of remarks}.
In section \ref{sec:linear_case}, we discuss the application of model reduction in the DtM framework   through the vehicle of a  two-dimensional model problem of the form \eqref{eq:model_problem}:
we present the HF FE discretization, and we introduce the lifting method to deal with Dirichlet BCs.
Then, in section \ref{sec:nonlinear_pb}, we apply our technique to the  incompressible Reynolds-averaged Navier-Stokes (RANS) simulations of  the flow past a two-dimensional  Ahmed  body   to demonstrate the applicability of the approach to nonlinear PDEs: as in the previous section, we present in detail the integration of the ROM with a HF stabilized FE model.
Finally, in section \ref{sec:conclusions}, we draw some conclusions and we outline a number of subjects of ongoing research.

\subsection{Relationship with previous works}

{
The idea of deforming the reduced mesh as opposed to restate the equations in the reference domains has been previously considered in several works.
Washabaugh et al.  \cite{washabaugh2016use}  have resorted to  
a similar discretize-then-map approach for steady aerodynamics flows in the finite volume framework: the authors 
resort to minimum residual projection and rely on Gappy POD, \cite{carlberg2013gnat}, for hyper-reduction.
Similarly, Dal Santo and Manzoni \cite{dal2019hyper} have resorted to a Galerkin ROM in combination with the (matrix) discrete empirical interpolation method (\cite{chaturantabut2010nonlinear,negri2015efficient}) for hyper-reduction.
However, in these works, the authors did not extensively compare the DtM approach with the MtD framework, and also did not discuss in detail the issue of dual residual  norm estimation, which is key to assess the accuracy of the ROM during the online stage.
}

EQ procedures are a class of hyper-reduction techniques, which aim to devise an efficient quadrature rule to approximate a suitable set of parameterized integrals associated with the reduced problem. 
We here rely on the EQ procedure in \cite{taddei2020space} that combines the approaches in \cite{yano2019discontinuous}  and  \cite{farhat2015structure}. 
As explained in Remark \ref{remark:RID} in the appendix, EQ procedures are related to the hyper-reduction procedure in \cite{ryckelynck2009hyper}: a thorough introduction to EQ procedures for hyper-reduction of projection-based ROMs is beyond the scope of the present work.
{
We refer to the recent reviews  \cite{farhat2021computational}
and  \cite{yano2021model} for a thorough discussion on EQ (also dubbed ``mesh sampling and weighting'); we further refer to 
\cite{chan2020entropy}  for  a discussion of the stability properties of the hyper-reduced system for hyperbolic  conservation laws.}

{
In section \ref{sec:nonlinear_pb}, we resort to 
 Least-Squares Petrov-Galerkin  (LSPG, \cite{carlberg2011efficient}) projection to devise an efficient ROM for the RANS equations.  As observed in \cite{grimberg2020stability,taddei2020space}, 
Petrov-Galerkin/minimum residual ROMs exhibit superior stability properties for non-coercive nonlinear PDEs. 
We remark that several authors have considered stabilized FE HF discretizations for ROM generation and construction, \cite{giere2015supg,pacciarini2014stabilized,stabile2019reduced}
: the above-mentioned papers, however, do not consider geometric parameterizations.
}

\section{Discretize-then-map treatment of parameterized geometries}
\label{sec:methods_DtM}
\subsection{Finite element spaces}
\label{sec:FEM_spaces}

For simplicity, 
we  omit dependence on parameter. 
Given the domain $\Omega \subset \mathbb{R}^D$,
we define the FE mesh $\mathcal{T} = \{ \texttt{D}_k  \}_{k=1}^{N_{\rm e}}$,  where $\texttt{D}_k \subset \Omega$ denotes  the $k$-th element of the mesh.
We define the reference element 
$\widehat{\texttt{D}} = \{ {X}\in [0,1]^D: \sum_{d=1}^D X_d<  1  \}$ 
 and the bijection 
${\Psi}_k$
from  $\widehat{\texttt{D}}$ to $\texttt{D}_k$ for $k=1,\ldots,N_{\rm e}$;  
we  denote the Jacobian and Jacobian determinant  of ${\Psi}_k$ by
${G}_k : = \nabla {\Psi}_k$ and $g_k: = {\rm det} ({G}_k)$.
Finally, we define the FE space of order \texttt{p} associated with the mesh $\mathcal{T}$:
\begin{equation}
\label{eq:FE_space}
\mathfrak{X}_{\mathcal{T}} : = \left\{
v\in C(\Omega) : \; \;
v    \circ {\Psi}_k \in \mathbb{P}_{\texttt{p}}(\widehat{\texttt{D}}),
\;
k=1,\ldots,N_{\rm e} 
\right\}.
\end{equation}

We consider a FE isoparametric discretization. We define the Lagrangian basis $\{ \ell_i \}_{i=1}^{n_{\rm lp}}$  of the polynomial space $\mathbb{P}_{\texttt{p}}(\widehat{\texttt{D}})$ associated with the nodes 
 $\{  {X}_i  \}_{i=1}^{n_{\rm lp}}$;
 then, we define the  mappings $\{ {\Psi}_k \}_k$ such that
  \begin{equation}
\label{eq:psi_mapping}
 {\Psi}_k( {{X}})
 =
 \sum_{i=1}^{n_{\rm lp}} \;
 {x}_{i,k}^{\rm hf}  \; \ell_i({X}),
\end{equation}
where $\{ {x}_{i,k}^{\rm hf}   := {\Psi}_k( {{X}}_i): \, i=1,\ldots,n_{\rm lp}, k=1,\ldots,N_{\rm e}  \}$ are the nodes of the mesh, and we define the basis functions
$\ell_{i,k} :=  \ell_i \circ  {\Psi}_k^{-1}: \texttt{D}_k \to \mathbb{R}$.  Note that ${\Psi}_k$ is completely characterized by  the nodes
in the $k$-th element 
 ${\texttt{X}}_{k}^{\rm hf} :=   \{  {x}_{i,k}^{\rm hf}  \}_{i=1}^{n_{\rm lp}}$, $k=1,\ldots,N_{\rm e}$.
We further introduce the
nodes of the mesh $\{ {x}_j^{\rm hf}  \}_{j=1}^{N_{\rm hf}}$ taken without repetitions  and the
 connectivity matrix $\texttt{T} \in \mathbb{N}^{n_{\rm lp}, N_{\rm e}}$,  such that
${x}_{i,k}^{\rm hf}   = {x}_{  \texttt{T}_{i,k} }^{\rm hf} $. 
{
We observe that we might use 
a different polynomial order for the mappings $\{ {\Psi}_k \}_k$: in view of mesh deformation (cf. \eqref{eq:psi_mapping_mapped}), it might be preferable for problems with straight boundaries to consider subparametric
discretizations based on linear elemental maps
(cf. section \ref{sec:optimal_convergence}).
}

Given the mesh $\mathcal{T}$ over $\Omega$ and the 
bijection ${\Phi}: \Omega \to {\Phi} ( \Omega )$, we  introduce the mapped mesh
 ${\Phi} (\mathcal{T})$ that shares with $\mathcal{T}$ the same connectivity matrix $\texttt{T}$ and has nodes $\{ {\Phi} ({x}_j^{\rm hf})    \}_{j=1}^{N_{\rm hf}}$.
 We denote by  $ {\Psi}_{k,\Phi}$ the elemental 
 mapping associated with the $k$-th element 
$\texttt{D}_{k,\Phi}$ 
 of  ${\Phi} (\mathcal{T})$; $ {\Psi}_{k,\Phi}$
  is given by
  \begin{equation}
\label{eq:psi_mapping_mapped}
{\Psi}_{k,\Phi}( {{X}})
 =
 \sum_{i=1}^{n_{\rm lp}} \;
 {\Phi}({x}_{i,k}^{\rm hf}  ) \; \ell_i({X}).
\end{equation}
  We also introduce the FE space
$\mathfrak{X}_{\Phi( \mathcal{T} )}$.

Given the vector $\mathbf{u} \in \mathbb{R}^{N_{\rm hf}}$, we define the corresponding FE field $u_{\Phi} \in \mathfrak{X}_{\Phi( \mathcal{T} )} $ such that $u_{\Phi} \in C \left( \bigcup_k \texttt{D}_{k,\Phi} \right)$ and
 \begin{subequations}
 \label{eq:vector2field}
  \begin{equation}
  u_{\Phi} \Big|_{  \texttt{D}_{k,\Phi}  }  \, = \,  
 \sum_{i=1}^{n_{\rm lp}} \; 
 \left( \mathbf{u}  \right)_{\texttt{T}_{i,k}} \; \,  \ell_{i,k,\Phi},
 \qquad
 k=1,\ldots,N_{\rm e},
 \end{equation}
 where $\ell_{i,k,\Phi} = \ell_i \circ  {\Psi}_{k,\Phi}^{-1}$.
If
${\Phi}$ is the identity map,   $\texttt{id}({x}) \equiv  {x}$, we use notation 
\begin{equation}
u = u_{\Phi = \texttt{id}} \in \mathfrak{X}_{\mathcal{T}}.
\end{equation}
  \end{subequations}
 Note that by changing ${\Phi}$ we effectively change the function field associated with the vector $\mathbf{u}$: we should thus interpret FE fields (and thus reduced-order spaces) as pairs of a FE vector and a mapping.

\subsection{MtD and DtM formulation of the Laplace problem}
In the DtM approach, we directly discretize 
\eqref{eq:model_problem}: given the mesh $\mathcal{T}$, the parameterized mapping ${\Phi}$ and the parameter $\mu\in \mathcal{P}$, we 
seek $u_{\mu,\Phi_{\mu}}^{\rm hf} \in \mathfrak{X}_{\mu}:=\mathfrak{X}_{\Phi_{\mu}(\mathcal{T})}$ such that
\begin{equation}
\label{eq:dtm_model_problem}
\int_{\Omega_{\mu}} \nabla    u_{\mu,\Phi_{\mu}}^{\rm hf}  \cdot \nabla v \; d{x} = 0 
\quad
\forall \, v \in \mathfrak{X}_{\mu} \cap   H_0^1( \Omega_{\mu}  ),
\quad
u_{\mu,\Phi_{\mu}}^{\rm hf} |_{\partial \Omega_{\mu}} = {h}_{\mu},
\end{equation}
with $\Omega_{\mu}  = {\Phi}_{\mu}(\Omega)$.
On the other hand, in the MtD approach, we discretize the mapped problem \eqref{eq:mapped_problem}: we seek 
$u_{\mu}^{\rm hf} \in \mathfrak{X} :=\mathfrak{X}_{\mathcal{T}}$ such that
\begin{equation}
\label{eq:mtd_model_problem}
\int_{\Omega} {K}_{\mu} \nabla u_{\mu}^{\rm hf} \cdot \nabla v \, d {x} \, = \,  0 \quad
\forall \, v \in \mathfrak{X} \cap H_0^1(\Omega),
\qquad
u_{\mu}^{\rm hf}|_{\partial \Omega} =\tilde{h}_{\mu}.
\end{equation}
In the remainder of this section, we compare the two formulations for the Laplace problem.

\subsubsection{Equivalence for piecewise-polynomial maps}
Exploiting the change-of-variable formula, we have
$$
\sum_{k=1}^{N_{\rm e}}
\;
\int_{\texttt{D}_k} \;
{K}_{\mu} \nabla w  \cdot \nabla v \, d {x} 
\;
=
\sum_{k=1}^{N_{\rm e}}
\;
\int_{ {\Phi}_{\mu}(\texttt{D}_k) } \;
\nabla w  \cdot \nabla v \, d {x}.
$$
Note   that ${\Psi}_{k,\Phi_{\mu}}$ \eqref{eq:psi_mapping_mapped} coincides with   ${\Phi} \circ {\Psi}_{k}$ in the nodes  
 $\{ {X}_i  \}_{i=1}^{n_{\rm lp}}$.  However, 
${\Psi}_{k,\Phi}$ and   ${\Phi}  \circ {\Psi}_{k}$  coincide in $\widehat{\texttt{D}}$ ---  
and thus $\texttt{D}_{k,\Phi} = {\Phi} (\texttt{D}_k)$ ---   only if 
 ${\Phi}  \circ {\Psi}_{k} \in \mathbb{P}_{\texttt{p}}$ for $k=1,\ldots,N_{\rm e}$. In conclusion,  we  obtain the  result contained in Lemma \ref{th:equivalence}: the proof   is a tedious but straightforward application of the chain rule and is here omitted.
 
\begin{lemma}
\label{th:equivalence}
Consider an isoparametric FE discretization of order \texttt{p}.
Then, \eqref{eq:dtm_model_problem} and \eqref{eq:mtd_model_problem} are equivalent if  
  ${\Phi}  \circ {\Psi}_{k} \in \mathbb{P}_{\texttt{p}}$ for $k=1,\ldots,N_{\rm e}$. 
  \end{lemma}
 
\subsubsection{Optimal convergence and discrete bijectivity}
\label{sec:optimal_convergence}
In order to highlight relevant features of the MtD and DtM approaches,
we consider the problem
\begin{subequations}
\label{eq:tricky_example}
\begin{equation}
-\partial_{x x} u =   \sin(\pi  x) \;\;\; x\in \Omega=(0,1),
\quad
u(0)=u(1)=0,
\end{equation}
which admits the unique solution $u(x) =  \frac{1}{\pi^2} \sin (\pi x)$. We further define the mapping
\begin{equation}
\Phi(x) =
\left\{
\begin{array}{ll}
\displaystyle{\frac{1}{2} x} & x \leq x_0, \\
\displaystyle{\frac{1}{2} \left( x_0 + 
\frac{2-x_0}{1-  x_0} (x - x_0) \right)  } & x \geq x_0 ; \\
\end{array}
\right.
\quad
x_0 =   \frac{1}{\sqrt{2}} .
\end{equation}
Clearly, $\Phi$ is a Lipschitz map from $\Omega$ in itself. We denote by $\mathcal{T}$ an uniform FE grid of $\Omega$ with $N_{\rm e}$ elements.  Problem \eqref{eq:tricky_example}  is intended to study the performance of the two strategies when combined with piecewise smooth maps, which are broadly used in MOR.
\end{subequations}

\begin{figure}[H]
\centering
\subfloat[]{
\begin{tikzpicture}[scale=0.7]
\begin{loglogaxis}[
xlabel = {\LARGE {$N_{\rm e}$} },
  ylabel = {\LARGE {$L^2$ error}},
  line width=1.2pt,
  mark size=3.0pt,
  ymin=0.000000001,   ymax=0.05,
 ]

\addplot[line width=1.pt,color=red,mark=square]  table {data/troubling_iso.dat}; \label{troubling:iso}
      
\addplot[line width=1.pt,color=blue,mark=triangle*] table {data/troubling_sub.dat};
\label{troubling:sub}
       
\addplot[line width=1.pt,color=black,mark=pentagon] table {data/troubling_mtd.dat};
\label{troubling:mtd}
   
\end{loglogaxis}
\end{tikzpicture}

\bgroup
\sbox0{\ref{data}}%
\pgfmathparse{\ht0/1ex}%
\xdef\refsize{\pgfmathresult ex}%
\egroup

}
~~
\subfloat[]{
\begin{tikzpicture}[scale=2.3]
\linethickness{0.3 mm}
\linethickness{0.3 mm}

\draw[ultra thick,dashed,color=blue]  (-1,0)--(0,0)--(1.2,-0.45);
\draw[thick,fill=gray, opacity=0.2]  (0,0)--(1,-0.375)--(-0.1,1)--(0,0);
\draw[ultra thick,fill=gray]  (0,0)--(1,-0.375)--(-0.6,0.3)--(0,0);
   \draw[color=black, fill=white] (-0.1,1) circle (.03);
\draw[color=black, fill=white] (-0.6,0.3) circle (.03);
\coordinate [label={above:  {\Huge {A}}}] (E) at (-0.1,1) ;
\coordinate [label={above:  {\Huge {B}}}] (E) at (-0.6,0.3) ;
\end{tikzpicture}
} 
\caption[Caption in ToC]{
optimal convergence and discrete bijectivity.
(a)
Performance of three FE discretizations for a Laplace problem.
Isoparametric  P3 DtM
\tikzref{troubling:iso};
subparametric  P1-P3 DtM
\tikzref{troubling:sub},
P3 MtD
\tikzref{troubling:mtd}.
(b)
Failure of the discretize-then-map approach for large deformations. The mapping maps the point A in B: this makes the resulting deformed mesh singular.
}
 \label{fig:troubling_case}
\end{figure}
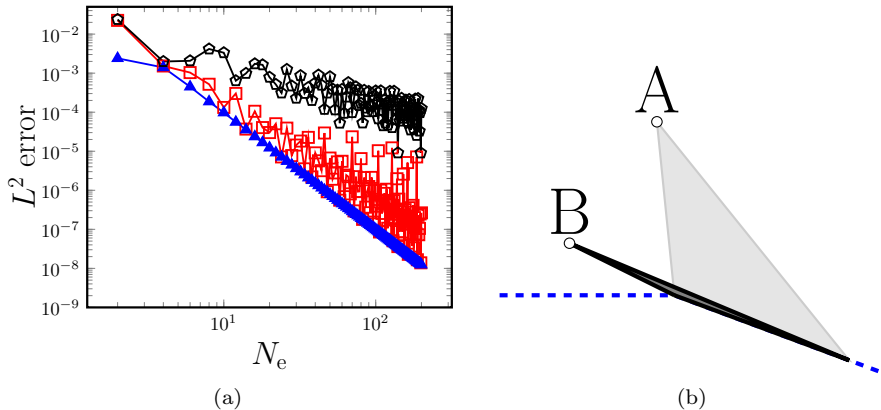

In Figure \ref{fig:troubling_case}(a), we show the behavior of the $L^2$ error associated with three FE discretizations:
a P3 isoparametric FE discretization of the   problem associated with $\Phi( \mathcal{T} )$ (\texttt{P3 DtM});
a P1-P3 subparametric   FE discretization of the   problem associated with $\Phi( \mathcal{T} )$ (\texttt{P1-P3 DtM} );
a P3 FE discretization of the mapped problem associated with $\mathcal{T}$ (\texttt{P3 MtD}).
We observe that the MtD approach fails to recover the optimal convergence rate ($r=4$) due to quadrature error and due to the fact that the MtD approach approximates the mapped field $u\circ \Phi$, which is considerably less smooth than $u$. The isoparametric discretization suffers due to stability issues: even if $\Phi$ is bijective, the local map $ \Psi_{k,\Phi} $ that contains $x_0$ is not necessarily bijective. On the other hand, the subparametric DtM approach recovers the optimal convergence rate as expected from standard FE theory (e.g., \cite[Proposition 3.4.1]{Quarteroni2008}). 

Note that for P1 discretizations in one-dimensional domains $\Phi(\mathcal{T})$ is guaranteed  to not have inverted elements; however,  the same result does not hold in higher dimensions.
Consider the element depicted in Figure \ref{fig:troubling_case}(b), the dashed line denotes the boundary of the domain: we can construct a mapping $\Phi$ that maps the point A in the point B in the Figure.
The deformed element thus reduces to a straight line regardless of the polynomial degree --- which implies that the corresponding elemental mapping $\Psi_{k,\Phi}$ is not bijective.
Note that for the MtD approach
 well-posedness depends on the fact that $K_{\mu}$ in \eqref{eq:mtd_model_problem}  is positive definite: the latter is independent of the underlying mesh and thus follows from the bijectivity of $\Phi$.
We conclude that, unlike MtD, the DtM approach might fail for smooth but large deformations, particularly in the presence of anisotropic elements.

\subsubsection{Application of hyper-reduction methods}
Hyper-reduction refers to a broad class of methods that aims to reduce the online assembling cost for projection-based ROMs; hyper-reduction techniques might be applied either at the 
continuous or at the discrete level --- \emph{continuous-based} and \emph{discrete-based} hyper-reduction. Continuous-based 
hyper-reduction methods (e.g., \cite{ballarin2016fast,ballarin2017numerical,grepl2007efficient}) aim 
to  devise  an affine approximation of the integrand  to ultimately obtain a parametrically-affine (cf. \cite{rozza2007reduced}) approximation of the problem:
$$
\mathcal{R}^{\rm aff}( u_{\mu}^{\rm hf}, v) 
=
\sum_{q=1}^{Q_{\rm a}} \; \Theta_{\mu}^q \mathcal{R}_q( u_{\mu}^{\rm hf},v) = 0,
$$
where $\mathcal{R}_1,\ldots,\mathcal{R}_{Q_{\rm a}}$ are parameter-independent forms that are polynomial in the first argument and linear in the second argument, and 
$\Theta_{\mu}^1,\ldots,\Theta_{\mu}^{Q_{\rm a}}$ are  parameter-dependent coefficients that can be 
rapidly evaluated. On the other hand, discrete-based hyper-reduction methods  (e.g., \cite{carlberg2013gnat,chaturantabut2010nonlinear,farhat2015structure,yano2019discontinuous}) 
aim to identify a ``reduced integration domain''. 

The MtD approach can cope with both continuous and discrete methods; furthermore, if the mapping ${\Phi}$ is piecewise-linear and discontinuities of ${\Phi}$ align with the elements' facets of the mesh, \ 
\eqref{eq:mtd_model_problem} admits a parametrically-affine decomposition without having to apply hyper-reduction.
On the other hand, the DtM approach can only cope with discrete-based hyper-reduction techniques, which require 
integration over a portion of the domain.

It is important to observe that parametrically-affine problems have been extensively studied in the model reduction literature and there exist specialized techniques to tackle them.
In particular, rigorous \emph{a posteriori} error estimation techniques are largely restricted to parametrically-affine problems and thus fit exclusively with the MtD framework.
\begin{itemize}
\item
Estimation of residual dual norm is straightforward for affine problems
(cf. \cite{rozza2007reduced}); furthermore, for certain classes of linear and nonlinear PDEs, the successive constraint method (SCM, \cite{huynh2007successive,huynh2010natural,yano2014space}) can be applied to obtain a rigorous error bound for the prediction error. In addition, the strategy in   \cite{schmidt2020rigorous} can later be applied to significantly sharpen the error bound.
\item
For a class of parametrically-affine PDEs, exact estimates can be derived to bound the error with respect to the exact solution 
(\cite{ali2017reduced,yano2016minimum}): this class of bounds is extremely important to devise spatio-parameter adaptive strategies for parametric problems,  \cite{yano2018reduced}.
 \end{itemize}
In this work, we adapt the approach in \cite{taddei2019offline} to devise a rapid and reliable estimator of the residual dual norm; furthermore, we show that the dual residual norm is highly correlated with the actual error: we can thus apply the approach in   \cite{drohmann2015romes}  to devise effective  probabilistic error bounds.

\subsection{Implementation for stabilized FE discretizations}
\label{sec:simplicity_implementation}
As stated in the introduction, practical implementation of the  MtD approach might be considerably more involved, particularly for stabilized FE formulations and nonlinear problems.
To motivate our claim, consider the streamline-upwind/ Petrov–Galerkin,
 (SUPG, \cite{brooks1982streamline})  discretization of the advection-diffusion equation
$$
-\Delta u + {b} \cdot \nabla u
= f \;\; {\rm in} \;\Omega_{\mu},
\quad
u|_{\partial \Omega_{\mu}} = 0;
$$
which can be stated as follows:
$$
\sum_{k=1}^{N_{\rm e}} 
\int_{\texttt{D}_{k,\Phi_{\mu}}}
\, 
\left(  \nabla u \cdot \nabla v \, + \,    {b} \cdot \nabla u
 v \, - \, f v  \right) \, d{x}
 \,+ \,  
\alpha h_{k,\Phi_{\mu}}
\int_{\texttt{D}_{k,\Phi_{\mu}}}
\, 
R(u) \frac{{b}}{\|{b} \|_2  } \cdot \nabla v     \, d{x}
= 0,
$$
where $R(u): = - \Delta u   \, + \,    {b} \cdot \nabla u - f$ is the strong residual, $\alpha>0$ is a user-defined constant, 
$h_{k,\Phi_{\mu}} = | \texttt{D}_{k,\Phi_{\mu}}  |^{1/D}$
and 
$ | \texttt{D}_{k,\Phi_{\mu}}  |$ is the volume of the element.

In order to apply the MtD approach, we should rewrite $- \Delta u$ in the reference configuration: if we denote by $\widetilde{\nabla}$ the gradient with respect to the reference coordinates, and we
define
$\tilde{u}:=u\circ {\Phi}_{\mu},
\tilde{v}:=u\circ {\Phi}_{\mu} $,
we find 
$$
\int_{\texttt{D}_{k,\Phi_{\mu}}}   \Delta u \, 
\left(  \frac{{b}}{\|{b} \|_2} \cdot \nabla v  \right)    \, d{x}
=
\int_{\texttt{D}_{k}} 
\; \left(
\left( \widetilde{\nabla} {\Phi}_{\mu} ^{-T} \;  \widetilde{\nabla} 
\right) \cdot 
\left( \widetilde{\nabla} {\Phi}_{\mu} ^{-T} \;  \widetilde{\nabla} \, \tilde{u}    
\right)
\right)
\frac{ \widetilde{\nabla}  {\Phi}_{\mu} ^{-1}  {b}}{\|{b} \|_2  } \cdot 
\left( \widetilde{\nabla}  \tilde{v}    \right)
\,
d {x}
 $$
Note that the expression at right-hand side requires the 
implementation of new local integration routines to deal with the mapped terms; in addition, it requires the  evaluation of the mapping ${\Phi}_{\mu}$, the gradient $ \widetilde{\nabla} {\Phi}_{\mu}$  and the Hessian 
 $ \widetilde{{H}} {\Phi}_{\mu}$ in the quadrature points. We envision that computation of second-order derivatives might be involved and  computationally expensive;  furthermore, second-order derivatives are not well-defined for Lipschitz maps. On the other hand, the DtM approach  only requires evaluation of ${\Phi}_{\mu}$ in the nodes of the mesh (cf. \eqref{eq:psi_mapping_mapped}); furthermore, it allows to reuse local integration routines that are used in the HF solver, cf. Algorithms \ref{alg:nonlinear_assembly} and \ref{alg:nonlinear_ROM_assembly} in section \ref{sec:nonlinear_pb}.

\section{Theory and results for the linear case}
\label{sec:linear_case}

We here present a  pMOR strategy  to estimate  elements of  the manifold 
$\mathfrak{M} :=\{\mathbf{u}_{\mu}^{\rm hf}: \mu \in \mathcal{P} \} \subset \mathfrak{U}$, where 
$\mathbf{u}_{\mu}^{\rm hf}$ denotes the  vector associated with the FE approximation to  \eqref{eq:model_problem} for a given $\mu \in \mathcal{P}$.
Given an estimate $\widehat{\mathbf{u}}_{\mu}$ of  $\mathbf{u}_{\mu}^{\rm hf}$ for some $\mu \in \mathcal{P} $, 
recalling the definitions  in \eqref{eq:vector2field},
we might define the FE fields 
$\widehat{u}_{\mu} \in  \mathfrak{X}_{\mathcal{T}}$ and $\widehat{u}_{\mu,\Phi_{\mu}} \in 
\mathfrak{X}_{\Phi_{\mu} ( \mathcal{T} )}
$: provided that the mapping is well-defined, we expect that
$\| {u}_{\mu}^{\rm hf} - \widehat{u}_{\mu,\Phi_{\mu}}  \|_{H^1(\Omega_{\mu})}
\approx 
\|\tilde{u}_{\mu}- \widehat{u}_{\mu}   \|_{H^1(\Omega)}$.

In section \ref{sec:preliminaries_linear}, we introduce useful notation and definitions.
In section \ref{sec:hf_model}, we present the HF formulation. 
Then, in section \ref{sec:galerkinROM}, we   derive a Galerkin ROM based on lifting: we discuss the offline/online computational decomposition, empirical interpolation of the Dirichlet BCs, the introduction of an  empirical quadrature 
rule { for the reduction of online costs to assemble the reduced matrix and vector}, and \emph{a posteriori} error estimation. 
In section  \ref{sec:model_problem},  we 
present the model problem considered in the numerical experiments, a two-dimensional potential flow past a parameterized  airfoil.
Finally,  in section \ref{sec:numerics}, we present  the 
numerical results  to prove the effectiveness of the approach.

\subsection{Preliminary definitions}
\label{sec:preliminaries_linear}
We omit dependence on parameter.
Given the Gaussian quadrature rule in $\widehat{\texttt{D}}$,
$\{  (\omega_q^{\rm g}, {X}_q^{\rm g}   )   \}_{q=1}^{n_{\rm q}}$, we define the high-fidelity quadrature (hfq)  rule in $\Omega$ as 
\begin{equation}
\label{eq:hf_quad}
\omega_{q,k}^{\rm hfq}:=\omega_q^{\rm g} \;\; g_k({X}_q^{\rm g}), \quad
{x}_{q,k}^{\rm hfq}  = {\Psi}_k ({X}_q^{\rm g}),
\qquad
q=1,\ldots,n_{\rm q}, \;\;
k=1,\ldots,N_{\rm e}.
\end{equation}
 Furthermore, we define $\mathbf{L}\in \mathbb{R}^{n_{\rm q}, n_{\rm lp}}$ and 
 $\mathbf{L}_{\nabla} \in \mathbb{R}^{n_{\rm q}, n_{\rm lp}, N_{\rm e}, D}$ such that
 \begin{equation}
 \label{eq:shape_functions}
 \left(  \mathbf{L} \right)_{q,i} = 
  \ell_i(  {X}_q^{\rm g} ),\quad
 \left(  \mathbf{L}_{\nabla} \right)_{q,i,k,d} = \partial_{x_d} \ell_{i,k}  (  {x}_{q,k}^{\rm hfq} ),
\end{equation}
with $q=1,\ldots,n_{\rm q}$, $i=1,\ldots,n_{\rm lp}$,
$k=1,\ldots,N_{\rm e}$,
$d=1,\ldots,D$.
Note that 
$ \ell_i(  {X}_q^{\rm g} ) = \ell_{i,k}({x}_{q,k}^{\rm hfq}  )$
for all $k=1,\ldots,N_{\rm e}$.
Quadrature rule and shape functions  are used for FE assembling; the HF  quadrature rule also provides  the baseline for the subsequent hyper-reduction procedure.
Given the mapping ${\Phi}$, we can introduce 
  the counterparts
 $\mathbf{L}_{\Phi}, 
\mathbf{L}_{\nabla, \Phi}$  of 
  $\mathbf{L}, 
\mathbf{L}_{\nabla}$ defined in \eqref{eq:shape_functions}: explicit expressions are obtained by replacing $\ell_{i,k} = \ell_i \circ {\Psi}_k^{-1}$ 
with  $\ell_{i,k,\Phi} = \ell_i \circ {\Psi}_{k,\Phi}^{-1}$ 
in \eqref{eq:shape_functions}: since
  $\mathbf{L}_{\Phi}= \mathbf{L}$ for any choice of ${\Phi}$, we shall omit the subscript $\Phi$ below.

We   define the Hilbert space $\mathfrak{U}  := ( \mathbb{R}^{N_{\rm hf}},  (\cdot,\cdot) )$ where the 
 inner product $(\cdot,\cdot)$ is given by 
\begin{equation}
\label{eq:inner_product_rom}
(\mathbf{w}, \mathbf{v}) = \int_{\Omega} \nabla w  \cdot \nabla v \, + \, 
w \cdot v  \; d {x}.
\end{equation}
We further define the induced norm
$\| \cdot \|  := \sqrt{(\cdot, \cdot)}$.
If we define the $H^1$ norm in $U \subset \mathbb{R}^D$ as 
$\|w \|_{H^1(U)} : = \sqrt{
\|w \|_{L^2(U)}^2 + \| \nabla w \|_{L^2(U)}^2}$, we obtain that
$\| \mathbf{u} \| = \|  u  \|_{H^1(\Omega)}$; furthermore, if ${\Phi}$ is a diffeomorphism, there exist $c,C>0$ such that
$
c \| \mathbf{u} \| \leq
\| u_{\Phi} \|_{H^1(\bigcup_k \texttt{D}_{k,\Phi}  )} \leq
C \| \mathbf{u} \|$ for all $\mathbf{u} \in \mathfrak{U}$.
   
We  define the ``elemental" restrictions
 such that
  $\mathbf{E}_k \mathbf{u} =[( \mathbf{u} )_{\texttt{T}_{1,k}},
\ldots,( \mathbf{u} )_{\texttt{T}_{n_{\rm lp},k}} ]^T$, 
for  $k=1, \ldots,  N_{\rm e}$: the restriction  $\mathbf{E}_k$ depends exclusively on the connectivity matrix $\texttt{T}$ and is thus independent of the mapping ${\Phi}$.
Given $\mathbf{u} \in \mathfrak{U}$, we define $\mathbf{u}^{\rm un} \in \mathbb{R}^{n_{\rm lp}, N_{\rm e}}$ such that
$\left(  \mathbf{u}^{\rm un}  \right)_{i,k} = 
\left(  \mathbf{E}_k \mathbf{u}   \right)_i$.
Similarly, given the reduced{-order} basis $\mathbf{Z} = [\boldsymbol{\zeta}_1,\ldots,\boldsymbol{\zeta}_N] \in \mathbb{R}^{N_{\rm hf}, N}$, we define the tensor
 $\mathbf{Z}^{\rm un} \in \mathbb{R}^{n_{\rm lp}, N_{\rm e}, N}$ such that
 \begin{equation}
 \label{eq:unassembled_reduced_basis}
\left( \mathbf{Z}^{\rm un}  \right)_{i,k,n}
=
\left( \boldsymbol{\zeta}_n^{\rm un}  \right)_{i,k}
=
\left( \mathbf{E}_k   \boldsymbol{\zeta}_n \right)_{i},
\quad
\left\{
\begin{array}{l}
i=1,\ldots,n_{\rm lp}, \\
k=1,\ldots,N_{\rm e},  \\
n=1,\ldots,N.\\
\end{array}
\right.
 \end{equation}

We denote by $\Gamma_{\rm dir} \subset \partial \Omega$ the Dirichlet boundary; in this work, we assume that each facet of the mesh is either strictly  contained in $\Gamma_{\rm dir}$ or has empty intersection with $\Gamma_{\rm dir}$. 
We introduce the set of indices of the mesh
$\texttt{I}_{\rm dir} \subset \{1, \ldots,N_{\rm hf}   \}$ belonging to  $\Gamma_{\rm dir} $; we denote by 
$M_{\rm hf} :=|  \texttt{I}_{\rm dir}  |$ the cardinality of the set. Then, 
we define the space
$\mathfrak{U}_0 :=\{ \mathbf{v} \in \mathfrak{U}: \mathbf{v} ( \texttt{I}_{\rm dir}) = \mathbf{0}  \}$ endowed with the norm $\|\cdot \|$, and the space
$\mathfrak{H} = (\mathbb{R}^{M_{\rm hf}}, \|  \cdot \|_{\rm dir} = \sqrt{(\cdot, \cdot)_{\rm dir}}   )$ where the norm $\|  \cdot \|_{\rm dir} $ is introduced below. Furthermore, we introduce the 
extension operator $\mathcal{H}: 
\mathfrak{H}  \to \mathfrak{U}$ such that for all $\mathbf{h} \in  \mathfrak{H}$
\begin{equation}
\label{eq:extension_operator}
\mathcal{H} \mathbf{h} (\texttt{I}_{\rm dir}) = \mathbf{h},
\qquad
 (  \mathcal{H} \mathbf{h},  \mathbf{v}  )
= 0
\quad
\forall \;  \mathbf{v} \in  \mathfrak{U}_0.
\end{equation}
$\mathcal{H}$ is used in section 
\ref{sec:galerkinROM} to derive the Galerkin ROM for the lifted solution field. Finally, we define the norm $\| \cdot \|_{\rm dir}$ such that for all $\mathbf{w} \in \mathfrak{U}$
\begin{equation}
\label{eq:bnd_norm}
\|  \mathbf{w} (\texttt{I}_{\rm dir}) \|_{\rm dir}
=
\sqrt{
\sum_{k=1}^{N_{\rm e}} \;
\int_{\partial \texttt{D}_k \cap \Gamma_{\rm dir}} \;
( w({x})  )^2 \, d {x}
}.
\end{equation}

\subsection{High-fidelity formulation of the model problem}
\label{sec:hf_model}

As in the previous section, we omit dependence on parameter.
 We discretize \eqref{eq:model_problem}
 using  a CG FE method: 
given the mesh $\mathcal{T}$,   the mapping ${\Phi}$, 
and the indices $\texttt{I}_{\rm dir} \subset \{1, \ldots,N_{\rm hf}   \}$ associated with $\Gamma_{\rm dir} =   \partial \Omega$, 
we denote by $\mathbf{h} \in \mathbb{R}^{M_{\rm hf}}$ the vector of boundary conditions
such that
$\left( \mathbf{h} \right)_i = h \left(
{\Phi} ( {x}_i^{\rm dir}   )
\right)$ with ${x}_i^{\rm dir} = {x}_{  ( \texttt{I}_{\rm dir}  )_i  }^{\rm hf}$, 
$i=1,\ldots, M_{\rm hf} :=|  \texttt{I}_{\rm dir}  |$.
Then,  
we seek $\mathbf{u}^{\rm hf} \in  \mathfrak{U}$  such that
$\mathbf{u}^{\rm hf}(\texttt{I}_{\rm dir}  ) = \mathbf{h}$ and 
\begin{subequations}
\label{eq:variational_formulation}
\begin{equation}
\mathcal{R}^{\rm hf}  (  \mathbf{u}^{\rm hf},  \mathbf{v} , {\Phi}  )
=
\sum_{k=1}^{N_{\rm e}} \;\;
r \left(
\mathbf{E}_k  \mathbf{u}^{\rm hf}, \;
\mathbf{E}_k \mathbf{v},  \; 
{\Phi} \left(   {\texttt{X}}_k^{\rm hf}     \right)    \right)
\; = \; 0
\quad
\forall \;  \mathbf{v} \in   \mathfrak{U}_0,
\end{equation}
where the local  residual $r$ is given by
\begin{equation}
\label{eq:mapped_residual_explicit_b}
\begin{array}{l}
\displaystyle{
r \left(
 \mathbf{w} , \mathbf{v} , {\Phi}
\left(   {\texttt{X}}_k^{\rm hf}     \right)    \right)
= 
\sum_{q=1}^{n_{\rm q}}
\; \omega_{q,k,\Phi}^{\rm hfq}  \, \cdot \, 
\left(
\nabla w_{\Phi}   \cdot \nabla v_{\Phi} 
\right)\big|_{ {x} = {x}_{q,k,\Phi}^{\rm hfq}    }
}\\[3mm]
=
\displaystyle{
\sum_{q=1}^{n_{\rm q}}
\; \omega_{q,k,\Phi}^{\rm hfq}  \, \cdot \, 
\sum_{j=1}^{n_{\rm lp}}
\left( \mathbf{v}   \right)_j \, \cdot \, 
\left(
\sum_{i=1}^{n_{\rm lp}}
\left( \mathbf{w}   \right)_i \, \cdot \, 
\left(
\sum_{d=1}^D
( \mathbf{L}_{\nabla, \Phi}   )_{q,i,k,d}  \, \cdot \, 
( \mathbf{L}_{\nabla, \Phi}   )_{q,j,k,d}  
\right)  
\right)
} \\[3mm]
=
\displaystyle{
\mathbf{v}^T 
\; \left(  \mathbf{A}_{k,\Phi}^{\rm un} \mathbf{w} \right)
},
\hfill
\forall \; 
 \mathbf{w} , \mathbf{v} \in \mathbb{R}^{n_{\rm lp}}.
 \\[3mm]
\end{array}
\end{equation}
Here,  $w_{{\Phi}},  v_{{\Phi}}$  are defined as in \eqref{eq:vector2field},  
and $ \mathbf{A}_{k, \Phi}^{\rm un} $  is given by
\begin{equation}
\label{eq:mapped_residual_explicit_c}
\left(  \mathbf{A}_{k, \Phi}^{\rm un}  \right)_{i,j}
=
\sum_{q=1}^{n_{\rm q}}
\; \omega_{q,k,\Phi}^{\rm hfq}  \, \cdot
\left(
\sum_{d=1}^D
( \mathbf{L}_{\nabla, \Phi}    )_{q,i,k,d}  \, \cdot \, 
( \mathbf{L}_{\nabla, \Phi}    )_{q,j,k,d}  
\right),
\;\;
i,j=1,\ldots,n_{\rm lp}.
\end{equation}
\end{subequations}
Given the solution $\mathbf{u}^{\rm hf}$ to \eqref{eq:variational_formulation}, we might define the corresponding FE fields  $u^{\rm hf}$ and $u_{\Phi}^{\rm hf}$ using \eqref{eq:vector2field}: the former is an approximation of the solution $\tilde{u}$ to  \eqref{eq:mapped_problem}; the latter is an approximation of the solution to 
\eqref{eq:model_problem}. 

Computation of $\left( \mathbf{L}_{\nabla, \Phi}  \right)_{\cdot, \cdot, k,\cdot} \in \mathbb{R}^{n_{\rm q}, n_{\rm lp}, D}$ for a given $k$
is dominated by the cost of computing    $\nabla {\Psi}_k$ in all quadrature points: the latter 
 involves the multiplication of a $D \times n_{\rm lp}$ matrix (associated with the nodes) with a $n_{\rm lp} \times n_{\rm q} D$ matrix (associated with the gradient of the shape functions). 
  Assuming $n_{\rm q} =\mathcal{O}( n_{\rm lp})$,
computation of  $\left( \mathbf{L}_{\nabla, \Phi}  \right)_{\cdot, \cdot, k,\cdot}$ 
 requires $\mathcal{O} ( n_{\rm lp}^2  )$ flops\footnote{We here assume that quadrature weights and quadrature points in $\widehat{\texttt{D}}$ and  shape functions  are computed in advance. 
Exploiting fast  matrix matrix multiplication algorithms, we might be able to assemble the tensor 
 $\mathbf{L}_{\nabla, \Phi}$ in 
 $\mathcal{O} ( n_{\rm lp}^{\beta} )$ flops
 with $2< \beta < 3$.
  }; similarly, computation of quadrature weights $\{  \omega_{q,k,\Phi}^{\rm hfq} \}_{q} $  and assembling of 
$\mathbf{A}_{k, \Phi}^{\rm un}$ 
using \eqref{eq:mapped_residual_explicit_c}
requires  
$\mathcal{O} ( n_{\rm lp}^3 )$ flops.
We remark that, thanks to a large body of research in the past few decades in the FE community,
efficient fully-vectorized routines are now available to compute local tensors of the form 
\eqref{eq:mapped_residual_explicit_c}
for structured and unstructured grids, for a broad range of PDEs; furthermore, the loop over the elements  to compute $\{  \mathbf{A}_{k,\Phi}^{\rm un} \}_k$ can be trivially parallelized.

\subsection{Galerkin reduced-order model for  the lifted solution field}
\label{sec:galerkinROM}
Given the reduced{-order} basis $\mathbf{Z} = [\boldsymbol{\zeta}_1,\ldots, \boldsymbol{\zeta}_N] \in \mathbb{R}^{N_{\rm hf}, N}$ such that
$\mathbf{Z}( \texttt{I}_{\rm dir}, :    ) = \mathbf{0}$, we aim to  approximate the solution $\mathbf{u}_{\mu}^{\rm hf}$ for a given $\mu\in \mathcal{P}$ as
\begin{equation}
\label{eq:galROM_abstract}
{\mathbf{u}}_{\mu}^{\star}
=
\mathbf{Z}   {\boldsymbol{\alpha}}_{\mu}^{\star}
+ \mathcal{H} \mathbf{h}_{\mu} ,
\quad
{\rm where} \;
\mathcal{R}_{\mu}^{\rm hf} ({\mathbf{u}}_{\mu}^{\star} , 
\boldsymbol{\zeta}_n; {\Phi}_{\mu}   )
= 0,
\;\;
n=1,\ldots,N.
\end{equation}
Note that computation of ${\mathbf{u}}_{\mu}^{\star}$ is highly inefficient: first, evaluation of 
$ \mathcal{H} \mathbf{h}_{\mu}$ requires the solution to a problem of size $N_{\rm hf}$; second, the assembling of the $N \times N$ system associated with \eqref{eq:galROM_abstract} requires integration over the whole FE mesh and thus scales with the total number of elements $N_{\rm e}$.

We shall resort to the  
empirical interpolation method (EIM, \cite{barrault2004empirical}) to obtain a low-dimensional   approximation of the Dirichlet datum of the form
\begin{equation}
\label{eq:affine_dirichlet}
\widehat{\mathbf{h}}_{\mu}
=
\mathbf{H} \; \mathbf{h}_{\mu} (\texttt{I}_{\rm ei}   ),
\;\;
{\rm where} \;\;
\mathbf{H}  \in \mathbb{R}^{  M_{\rm hf}, M }, \;\; 
\texttt{I}_{\rm ei}  \subset \{1,\ldots,M_{\rm hf} \},
\;\;
| \texttt{I}_{\rm ei}  | = M \ll M_{\rm hf}.
\end{equation}
Since $\mathcal{H}$ is parameter-independent, we can precompute the surrogate basis 
$\mathbf{W} = [\boldsymbol{\psi}_1,\ldots,\boldsymbol{\psi}_M] 
= \mathcal{H} \mathbf{H}
\in \mathbb{R}^{N_{\rm hf}, M}$; finally, we define the approximate lift 
$\widehat{\mathbf{e}}_{\mu} = \mathbf{W}  \mathbf{h}_{\mu} (\texttt{I}_{\rm ei}   )$.

Exploiting \eqref{eq:mapped_residual_explicit_c}, it is easy to verify that 
\begin{subequations}
\label{eq:galROM_hf}
\begin{equation}
\mathcal{R}_{\mu}^{\rm hf} \left(
\mathbf{Z}  \widehat{\boldsymbol{\alpha}}_{\mu}^{\rm hf} + \widehat{\mathbf{e}}_{\mu} , 
\boldsymbol{\zeta}_n; {\Phi}_{\mu}   \right)
= 0 
\; \Leftrightarrow \;
\widehat{\mathbf{A}}_{\mu}^{\rm hf} 
\, \widehat{\boldsymbol{\alpha}}_{\mu}^{\rm hf} 
\, = \,
\widehat{\mathbf{F}}_{\mu}^{\rm hf} 
\end{equation}
where
\begin{equation}
\label{eq:galROM_hf_b}
\left\{
\begin{array}{ll}
\displaystyle{
\widehat{\mathbf{A}}_{\mu}^{\rm hf}  :=
\sum_{k=1}^{N_{\rm e}} \;
\widehat{\mathbf{A}}_{k, \mu} },
&
\displaystyle{
\left( \widehat{\mathbf{A}}_{k, \mu} \right)_{n,n'} =
 \left( \mathbf{E}_k  \boldsymbol{\zeta}_{n} \right)^T
\mathbf{A}_{k, \Phi_{\mu} }^{\rm un} \, 
 \mathbf{E}_k  \boldsymbol{\zeta}_{n'}  
 }
\\[4mm]
\displaystyle{
\widehat{\mathbf{F}}_{\mu}^{\rm hf}  :=
\sum_{k=1}^{N_{\rm e}} \;
\widehat{\mathbf{F}}_{k, \mu}},
&
\displaystyle{
\left( \widehat{\mathbf{F}}_{k, \mu} \right)_{n} =
-  \left( \mathbf{E}_k    \boldsymbol{\zeta}_{n} \right)^T
\mathbf{A}_{k, \Phi_{\mu} }^{\rm un} \, 
 \mathbf{E}_k   \widehat{\mathbf{e}}_{\mu} 
}.
\\
\end{array}
\right.
\end{equation}
\end{subequations}
Assuming that $\{  \mathbf{E}_k \mathbf{Z} \}_k$ and 
$\{  \mathbf{E}_k \mathbf{W} \}_k$ are precomputed during the offline stage, the cost of assembling the reduced system scales with 
$\mathcal{O} ( 
N_{\rm e}  n_{\rm lp} ( ( N+n_{\rm lp}    )^2 + M )
 )$. To reduce assembling costs, we replace the HF residual with the weighted  residual $\mathcal{R}_{\mu}^{\rm eq}$, 
\begin{equation}
\label{eq:weighted_residual}
\mathcal{R}_{\mu}^{\rm eq}  (  \mathbf{w},  \mathbf{v} , {\Phi}_{\mu}  )
=
\sum_{k=1}^{N_{\rm e}} \;\; \rho_k^{\rm eq}
r \left(
\mathbf{E}_k  \mathbf{w} , \; \mathbf{E}_k \mathbf{v}, {\Phi}_{\mu}
\left(   {\texttt{X}}_k^{\rm hf}     \right)    \right),
\end{equation}
where $\boldsymbol{\rho}^{\rm eq} = [\rho_1^{\rm eq}, \ldots,\rho_{N_{\rm eq}}^{\rm eq} ]^T$  is a sparse vector of positive weights such that $\|  \boldsymbol{\rho}^{\rm eq}  \|_{\ell^0} = Q \ll N_{\rm e}$. We denote by $\texttt{I}_{\rm eq} = \{ k \in \{ 1,\ldots,N_{\rm e} \} : \rho_k^{\rm eq}>0   \}$ the set of indices associated with the sampled elements.

In conclusion, 
given $\mu \in \mathcal{P}$, we define the estimate 
$\widehat{\mathbf{u}}_{\mu}
=
\mathbf{Z} \widehat{\boldsymbol{\alpha}}_{\mu} +
\widehat{\mathbf{e}}_{\mu}
$ of 
the solution $\mathbf{u}_{\mu}^{\rm hf}$ to  
\eqref{eq:variational_formulation},
such that
\begin{subequations}
\label{eq:GalROM}
\begin{equation}
\label{eq:GalROMa}
\widehat{\mathbf{e}}_{\mu} = \mathbf{W}  \; \mathbf{h}_{\mu} (\texttt{I}_{\rm ei}   ),
\qquad
\mathcal{R}_{\mu}^{\rm eq}  ( 
\mathbf{Z} \widehat{\boldsymbol{\alpha}}_{\mu} +
\widehat{\mathbf{e}}_{\mu},  \boldsymbol{\zeta}_n , {\Phi}_{\mu}  ) = 0,
\;\;
n=1,\ldots,N.
\end{equation}
Exploiting previous identities, we find that $\widehat{\boldsymbol{\alpha}}_{\mu} $ satisfies
\begin{equation}
\label{eq:GalROMb}
\widehat{\mathbf{A}}_{\mu}^{\rm eq} 
\, \widehat{\boldsymbol{\alpha}}_{\mu} 
\, = \,
\widehat{\mathbf{F}}_{\mu}^{\rm eq},
\qquad
{\rm where} \;\; 
\widehat{\mathbf{A}}_{\mu}^{\rm eq}  =
\sum_{k \in \texttt{I}_{\rm eq}  }  \; 
\rho_k^{\rm eq} \; 
\widehat{\mathbf{A}}_{k, \mu},
\quad
\widehat{\mathbf{F}}_{\mu}^{\rm eq}  =
\sum_{k \in \texttt{I}_{\rm eq}  }  \; \rho_k^{\rm eq} \; 
\widehat{\mathbf{F}}_{k, \mu}.
\end{equation}
\end{subequations}
Online assembling of the reduced system requires the computation of the $Q \cdot n_{\rm lp}$ nodes of the mapped mesh,
$\{ {\Phi}_{\mu}({x}_{ i,k}^{\rm hf}): i=1,\ldots,n_{\rm lp}, k \in  \texttt{I}_{\rm eq} \}$,  
the assembling of the elemental matrices
$\{ \mathbf{A}_{k, \Phi_{\mu} }^{\rm un}   \}_{k \in   \texttt{I}_{\rm eq}}$ in \eqref{eq:mapped_residual_explicit_c},
the 
computation of $\{ \widehat{\mathbf{A}}_{k, \mu}   \}_k$ in \eqref{eq:galROM_hf_b}, 
 and then the summation over the sampled elements. Provided that evaluation of ${\Phi}_{\mu}$ in a given point requires $\mathcal{O}(1)$ flops, the overall assembling cost scales with $\mathcal{O} ( 
Q  n_{\rm lp} ( ( N+n_{\rm lp}    )^2 + M )
 )$.

\begin{algorithm}[H]                      
\caption{Offline-online computational decomposition}     
\label{alg:off_on}     

\textbf{Offline stage}
\medskip

\begin{algorithmic}[1]

\State
Compute samples of the HF solution manifold 
$\{ \mathbf{u}_{\mu^k}^{\rm hf}   \}_{k=1}^{n_{\rm train}}$ and compute the lifted fields
$\{  \mathring{\mathbf{u}}_{\mu^k}^{\rm hf} = 
\mathbf{u}_{\mu^k}^{\rm hf}  - \mathcal{H} 
{\mathbf{h}}_{\mu^k}    \}_k$

\State
Compute 
the POD space
$\mathbf{Z} = \texttt{POD} \left(
\{ \mathring{\mathbf{u}}_{\mu^k}^{\rm hf}  \}_k, (\cdot, \cdot), tol_{\rm pod} \right)$ .
\medskip

\State
Compute $\mathbf{H}$ and $\texttt{I}_{\rm ei}$ 
in \eqref{eq:affine_dirichlet}  using EIM
and then set 
$\mathbf{W} = \mathcal{H} \mathbf{H}$ 
 (cf. section \ref{sec:EIM}).
\medskip

\State
Find the sparse vector of weights $\boldsymbol{\rho}^{\rm eq}$ in \eqref{eq:weighted_residual}
(cf. section \ref{sec:EQP}).
\medskip

\State
Construct the offline structures for dual residual norm  calculation (cf.  section  \ref{sec:a_posteriori}).
\end{algorithmic}
\medskip

\textbf{Online stage} (for any given $\mu \in \mathcal{P}$)
\medskip

\begin{algorithmic}[1]
\State
{ Compute the deformed nodes
$\{ {\Phi}_{\mu}({x}_{ i,k}^{\rm hf}): i=1,\ldots,n_{\rm lp}, k \in  \texttt{I}_{\rm eq} \}$}.
\medskip

\State
Assemble
$\widehat{\mathbf{A}}_{\mu}^{\rm eq}$ and
$\widehat{\mathbf{F}}_{\mu}^{\rm eq}$ in \eqref{eq:GalROMb}.
\medskip

\State
Solve $\widehat{\mathbf{A}}_{\mu}^{\rm eq} 
\, \widehat{\boldsymbol{\alpha}}_{\mu} 
\, = \,
\widehat{\mathbf{F}}_{\mu}^{\rm eq}$ and return
$
 \widehat{\mathbf{u}}_{\mu}=
\mathbf{Z} \widehat{\boldsymbol{\alpha}}_{\mu} + \mathbf{W}  \mathbf{h}_{\mu} (\texttt{I}_{\rm ei}   )$.
\medskip

\State
Estimate the dual residual  norm \eqref{eq:dual_residual}  (cf. section   \ref{sec:a_posteriori}).
\end{algorithmic}

\end{algorithm}

In Algorithm \ref{alg:off_on}, we outline the offline/online computational decomposition considered in the numerical experiments; the procedure includes the offline construction and the online evaluation of the dual residual estimator that is used to assess the accuracy of the ROM. 
Computation of \eqref{eq:GalROM} requires the storage of the reduced{-order} bases $\mathbf{Z}, \mathbf{W}$ and of 
the nodes $\{ {x}_{i,k}^{\rm hf} : i=1,\ldots,n_{\rm lp}, k \in \texttt{I}_{\rm eq}   \}$ in the selected elements 
($n_{\rm lp} Q (M+N+D)$ doubles),  and the storage of quantities 
(shape functions, quadrature rule)
associated with the master element 
$\widehat{\texttt{D}}$ ($\mathcal{O}(n_{\rm lp} n_{\rm q})$ doubles )
--- computational and storage costs do not consider dual residual norm estimation, which is addressed in section \ref{sec:a_posteriori}.

{ We here resort to  a non-adaptive (deterministic or random) sampling of the parameter domain and we resort to proper orthogonal decomposition 
(POD, \cite{volkwein2011model})  to generate the reduced{-order}  basis for the ROM.
The size of the spaces is chosen according to the criterion
\begin{equation}
\label{eq:POD_cardinality_selection}
N := \min \left\{
N': \, \sum_{n=1}^{N'} \lambda_n \geq  \left(1 - tol_{\rm pod} \right) 
\sum_{i=1}^{n_{\rm train}} \lambda_i
\right\},
\end{equation} 
where $\{  \lambda_i  \}_i$ are the eigenvalues in descending order of the POD covariance matrix, and 
$tol_{\rm pod} $ is a suitable tolerance.
Given snapshots $\{  \mathbf{w}^k \}_k \subset \mathbb{R}^{N_{\rm pod}}$, 
 notation $\mathbf{W} = \texttt{POD} \left(
\{  \mathbf{w}^k \}_k, (\cdot, \cdot)_{\rm pod}, tol_{\rm pod} \right)$ signifies that 
$\mathbf{W} \in \mathbb{R}^{N_{\rm pod} , N}$ is constructed using POD  based on the inner product 
$(\cdot, \cdot)_{\rm pod}$ and tolerance 
$tol_{\rm pod} $ in \eqref{eq:POD_cardinality_selection}.
}

In the next three sections, we review how to construct the approximate lift $\widehat{\mathbf{e}}$,  
the weights $\boldsymbol{\rho}^{\rm eq} \in \mathbb{R}^{N_{\rm e}}$ in  
\eqref{eq:weighted_residual}, 
and how to estimate the dual residual error indicator;
 we state upfront that a thorough discussion and analysis of EIM and EQ procedures is beyond the scope of the present work: we refer to 
\cite{barrault2004empirical} (see also \cite{chaturantabut2010nonlinear})
and to the reviews  \cite{hesthaven2016certified,quarteroni2015reduced}
for a detailed presentation of EIM; on the other hand, we refer to the previously mentioned papers  \cite{yano2019discontinuous,farhat2015structure,farhat2021computational,yano2021model}
for  a discussion on EQ procedures for linear and nonlinear problems.

\subsubsection{Affine approximation of the Dirichlet datum}
\label{sec:EIM}

Given the parameters $\{  \mu^k \}_{k=1}^{n_{\rm train}} \subset \mathcal{P}$, we compute
$\{ \mathbf{h}^k = \mathbf{h}_{\mu^k}   \}_k$ and we apply POD based on the $(\cdot, \cdot)_{\rm dir}$ inner product to obtain a $M$-dimensional approximation space 
$\mathfrak{H}_M 
= {\rm span} \{  \boldsymbol{\xi}_m \}_{m=1}^M
\subset 
\mathfrak{H}_{\rm dir}$ 
where the integer $M$ is chosen according to \eqref{eq:POD_cardinality_selection}.  
Then, we pass $\boldsymbol{\Xi} :=[ \boldsymbol{\xi}_1,\ldots, \boldsymbol{\xi}_M ]$ to the EIM Greedy algorithm
(see, e.g., \cite[Algorithm 10.3]{quarteroni2015reduced}) to find the set of indices $\texttt{I}_{\rm ei} \subset \{1,\ldots,M_{\rm hf} \}$, $| \texttt{I}_{\rm ei}  | = M$. Finally, we define the matrix $\mathbf{W}$ in  \eqref{eq:GalROMa} as
$\mathbf{W} = \mathcal{H}  \boldsymbol{\Xi} \,  ( \boldsymbol{\Xi}(\texttt{I}_{\rm ei} ,:) ) ^{-1}$.

{
We might avoid the use of EIM to approximate the Dirichlet datum by resorting to a weak imposition of the Dirichlet condition ---
see, e.g., \cite{yano2019discontinuous} and also 
 \cite{babuvska1973finite,juntunen2009nitsche} for a discussion in the continuous setting.
In this case, we can use EQ to simultaneously deal with hyper-reduction and BC approximation. Since we do not pursue this strategy in this paper, we do not further comment on this issue.
}

\subsubsection{Empirical quadrature procedure}
\label{sec:EQP}
Following \cite{yano2019discontinuous}, we seek  $\boldsymbol{\rho}^{\rm eq} \in \mathbb{R}_+^{N_{\rm e}}$ such that
(i) the number of nonzero entries in 
 $\boldsymbol{\rho}^{\rm eq}$,  $\| \boldsymbol{\rho}^{\rm eq}   \|_{\ell^0}$,   is as small as possible;
(ii, \emph{constant function constraint})
 the constant function is approximated correctly in $\Omega$ (i.e., ${\Phi} = \texttt{id}$), 
 \begin{equation}
 \label{eq:constant_function_constraint}
\Big|
\sum_{k=1}^{N_{\rm e}} \rho_k^{\rm eq} | \texttt{D}_k  |
\,-\,
| \Omega | 
 \Big|
 \ll 1; 
 \end{equation}
(iii, \emph{manifold accuracy constraint})
for all $\mu \in \mathcal{P}_{\rm train} = \{  \mu^k \}_{k=1}^{n_{\rm train}}$, the  empirical residual satisfies
\begin{equation}
\label{eq:accuracy_constraint}
\Big \|
\left(  \widehat{\mathbf{A}}_{\mu}^{\rm hf}   \right)^{-1}
\left(  
\widehat{\mathbf{A}}_{\mu}^{\rm eq}  \, \widehat{\boldsymbol{\alpha}}_{\mu} ^{\rm hf} 
- \widehat{\mathbf{F}}_{\mu}^{\rm eq} 
\right)
 \Big \|_2
 \ll 1.
\end{equation}
We refer to  \cite{yano2019discontinuous}  for a detailed justification of \eqref{eq:constant_function_constraint} and \eqref{eq:accuracy_constraint} for linear and nonlinear PDEs. 

It is easy to verify that \eqref{eq:constant_function_constraint} and \eqref{eq:accuracy_constraint} could be rewritten in matrix form as
$$
\Big|
\sum_{k=1}^{N_{\rm e}} \rho_k^{\rm eq} | \texttt{D}_k  |
\,-\,
| \Omega | 
\Big|
=
|  \mathbf{G}_{\rm const}  \boldsymbol{\rho}^{\rm eq}  - | \Omega |  
\Big|
 \ll 1,
 $$
 with $\mathbf{G}_{\rm const} = [  | \texttt{D}_1  |, \ldots,  | \texttt{D}_{N_{\rm e}}  |  ]$, and 
 $$
 \Big \|
\left(  \widehat{\mathbf{A}}_{\mu}^{\rm hf}   \right)^{-1}
\left(  
\widehat{\mathbf{A}}_{\mu}^{\rm eq}  \, \widehat{\boldsymbol{\alpha}}_{\mu} ^{\rm hf} 
- \widehat{\mathbf{F}}_{\mu}^{\rm eq} 
\right)
 \Big \|_2
 =
  \Big \|
 \mathbf{G}_{\mu}   \, \boldsymbol{\rho}^{\rm eq}   \Big \|_2
 \ll 1,
 $$
 with
 $
{\mathbf{G}}_{\mu}  
=
  \left(  \widehat{\mathbf{A}}_{\mu}^{\rm hf}   \right)^{-1}
  \left[
    \left( 
  \widehat{\mathbf{A}}_{1, \mu}
  \widehat{\boldsymbol{\alpha}}_{\mu} ^{\rm hf} 
  -\widehat{\mathbf{F}}_{1, \mu}   \right),\ldots,
    \left( 
  \widehat{\mathbf{A}}_{N_{\rm e}, \mu}
  \widehat{\boldsymbol{\alpha}}_{\mu} ^{\rm hf} 
  -\widehat{\mathbf{F}}_{N_{\rm e}, \mu}   \right)
  \right]
 $.
The problem of finding   $\boldsymbol{\rho}^{\rm eq}$  can thus be reformulated as a sparse-representation (or best-subset selection) problem:
\begin{subequations}
\label{eq:sparse_representation}
\begin{equation}
\min_{  \boldsymbol{\rho} \in \mathbb{R}^{N_{\rm e}} }
\;
\| \boldsymbol{\rho}   \|_{\ell^0},
\quad
{\rm s.t} \quad
\left\{
\begin{array}{l}
\|\mathbf{G} \boldsymbol{\rho} - \mathbf{b}  \|_{\star} \leq \delta; \\[3mm]
\boldsymbol{\rho} \geq \mathbf{0}; \\
\end{array}
\right.
\end{equation}
for   suitable choices of the   vector norm 
$\|\cdot   \|_{\star} $ and the tolerance $\delta>0$, and
\begin{equation}
\mathbf{G}  = \left[
\begin{array}{l}
\mathbf{G}_{\rm const}   \\[3mm]
\mathbf{G}_{\mu^1}  \\
\vdots \\
\mathbf{G}_{\mu^{n_{\rm train}}}  \\
\end{array}
\right],
\;\;
\mathbf{b}  = \left[
\begin{array}{l}
|\Omega|  \\[3mm]
\mathbf{0}  \\
\end{array}
\right]
\end{equation}
\end{subequations}
A large body of research focuses on the development of robust algorithms to find approximate 
solutions to \eqref{eq:sparse_representation}.
Following \cite{farhat2015structure}, we resort  to the nonnegative linear least-squares method (see \cite{lawson1974solving}) 
---  more in detail, we rely on 
the Matlab function \texttt{lsqnonneg},
which takes as input the pair $(\mathbf{G}, \mathbf{b})$ and a tolerance $tol_{\rm eq}>0$ and returns  the sparse vector $\boldsymbol{\rho}^{\rm eq}$,
$$
[\boldsymbol{\rho}^{\rm eq}] =  \texttt{lsqnonneg} \left( \mathbf{G}, \mathbf{b}, tol_{\rm eq} \right).
$$
We remark that for large-scale problems  computation of the solution to the  nonnegative linear least-squares problem might be prohibitively expensive: to address this issue, efficient partitioned approaches have been developed in \cite{grimberg2020mesh}.

\subsubsection{A posteriori error estimation}
\label{sec:a_posteriori}

Due to the absence of reliable a priori error estimates, it is of paramount importance  to assess the accuracy of the ROM during the online stage; for this reason, a large body of work has been devoted to the development of rapid and reliable \emph{a posteriori} error indicators { --- see  \cite{rozza2007reduced,hesthaven2016certified,quarteroni2015reduced}  and the references therein}. Most proposals rely on the computation (estimation) of the dual residual norm
\begin{equation}
\label{eq:dual_residual}
\mathfrak{R}_{\mu}^{\rm hf} ( \widehat{\mathbf{u}}_{\mu}    )
:=
\sup_{ \mathbf{v} \in \mathfrak{U}_0 } \;
\frac{\mathcal{R}_{\mu}^{\rm hf}( \widehat{\mathbf{u}}_{\mu}, \mathbf{v}; {\Phi}_{\mu} )}{\| \mathbf{v} \| }.
\end{equation}
Here, we propose to adapt the strategy in \cite{taddei2019offline} to devise a rapid and reliable estimator of the dual residual norm \eqref{eq:dual_residual}: the approach involves the definition of a low-dimensional  \emph{empirical test space} and of an empirical  quadrature rule.

We approximate the dual residual norm
$\mathfrak{R}_{\mu}^{\rm hf} ( \widehat{\mathbf{u}}_{\mu}  )$ for any given $\mu$ as 
\begin{equation}
\label{eq:dual_norm_estimate}
\widehat{\mathfrak{R}}_{\mu}^{\rm eq}
:=
\sup_{ \mathbf{v} \in \mathfrak{Y}_{J_{\rm r}}^{\rm r} } \;
\frac{\mathcal{R}_{\mu}^{\rm eq,r}( \widehat{\mathbf{u}}_{\mu}, \mathbf{v}; {\Phi}_{\mu} )}{\| \mathbf{v} \| },
\end{equation}
where $\mathfrak{Y}_{J_{\rm r}}^{\rm r} = {\rm span} \{   \boldsymbol{\eta}_j \}_{j=1}^{J_{\rm r}} \subset \mathfrak{U}_0$ is referred to as \emph{empirical test space} and 
$\mathcal{R}_{\mu}^{\rm eq,r}$ satisfies
\begin{equation}
\label{eq:weighted_residual_res}
\mathcal{R}_{\mu}^{\rm eq,r} (  \mathbf{w},  \mathbf{v} , {\Phi}_{\mu}  )
=
\sum_{k=1}^{N_{\rm e}} \;\; \rho_k^{\rm eq,r}
r \left(
\mathbf{E}_k  \mathbf{w} , \; \mathbf{E}_k \mathbf{v}, {\Phi}_{\mu}
\left(   {\texttt{X}}_k^{\rm hf}     \right)    \right),
\quad
\forall \, \mathbf{w}, \mathbf{v} \in \mathfrak{U}, \;\; 
\mu\in \mathcal{P}.
\end{equation}
Similarly to \eqref{eq:weighted_residual},
$\boldsymbol{\rho}^{\rm eq,r} = [\rho_1^{\rm eq}, \ldots,\rho_{N_{\rm eq}}^{\rm eq} ]^T$  is a sparse vector of positive weights such that $\|  \boldsymbol{\rho}^{\rm eq,r}  \|_{\ell^0} = Q_{\rm r} \ll N_{\rm e}$;  we denote by $\texttt{I}_{\rm eq,r} = \{ k \in \{ 1,\ldots,N_{\rm e} \} : \rho_k^{\rm eq,r}>0   \}$ the set of indices associated with the sampled elements. Provided that $\{   \boldsymbol{\eta}_j \}_{j=1}^{J_{\rm r}}$ is an orthonormal basis of $ \mathfrak{Y}_{J_{\rm r}}^{\rm r}$, we find that
\begin{subequations}
\label{eq:residual_calculation}
\begin{equation}
\label{eq:residual_calculation_a}
\widehat{\mathfrak{R}}_{\mu}^{\rm eq} = \| \widehat{\mathbf{R}}_{\mu}^{\rm eq,r}  \|_2 ,
\quad
\widehat{\mathbf{R}}_{\mu}^{\rm eq,r} = \sum_{k \in \texttt{I}_{\rm eq,r}   }
\, \rho_k^{\rm eq,r} \; \widehat{\mathbf{R}}_{k,\mu}^{\rm r},
\end{equation}
where $\{  \widehat{\mathbf{R}}_{k,\mu}^{\rm r}  \}_k \subset \mathbb{R}^{J_{\rm r}}$ are given by
\begin{equation}
\label{eq:residual_calculation_b}
\left( \widehat{\mathbf{R}}_{k,\mu}^{\rm r} \right)_j:=
r \left( \mathbf{E}_k \widehat{\mathbf{u}}_{\mu}    , \mathbf{E}_k  \boldsymbol{\eta}_j , {\Phi}_{\mu}
\left(   {\texttt{X}}_k^{\rm hf}     \right)    \right)
=
 \left(   \left( \mathbf{E}_k   \boldsymbol{\eta}_j \right)^T
\mathbf{A}_{k, \Phi_{\mu} }^{\rm un} \, 
\mathbf{E}_k \, 
\left(      \mathbf{Z} \widehat{\boldsymbol{\alpha}}_{\mu} +  \widehat{\mathbf{e}}_{\mu}   \right)      \right)_{j},
\end{equation}
for $j=1,\ldots,J_{\rm r}$ and  $k=1,\ldots,N_{\rm e}$.
\end{subequations}
Given $\mu\in \mathcal{P}$ and the estimate $\widehat{\boldsymbol{\alpha}}_{\mu}$, 
the cost of estimating $\widehat{\mathfrak{R}}_{\mu}^{\rm eq}$ is dominated 
by the cost of assembling the local residuals  $\{  \widehat{\mathbf{R}}_{k,\mu}^{\rm r}  \}_k$ and scales with\footnote{
Computation of $\{  \mathbf{A}_{k,\Phi_{\mu}} \}_{k}$:
$\mathcal{O} (  n_{\rm lp}^3 Q_{\rm r} )$;
computation of $\mathbf{E}_k \widehat{\mathbf{u}}_{\mu}$
(given $\widehat{\boldsymbol{\alpha}}_{\mu}$):
$\mathcal{O} (  n_{\rm lp}  Q_{\rm r} (N+M) )$;
computation of $\widehat{\mathbf{R}}_{k,\mu}^{\rm r}$:
$\mathcal{O} (  Q_{\rm r}  ( n_{\rm lp}^2  + J_{\rm r} n_{\rm lp} )  )$;
computation of  $\| \widehat{\mathbf{R}}_{\mu}^{\rm r}\|_2$:
$\mathcal{O} (     ( Q_{\rm r} +1)  J_{\rm r}  )$.
} 
$\mathcal{O} \left(  
Q_{\rm r} \, \left(  
n_{\rm lp}^3 + n_{\rm lp} (N+M) +   J_{\rm r} n_{\rm lp} \right)
 \right)$. In the remainder of this section, we illustrate how to construct the test space $\mathfrak{Y}_{J_{\rm r}}^{\rm r}$ and the quadrature weights $\boldsymbol{\rho}^{\rm eq,r} $.

We shall choose $\mathfrak{Y}_{J_{\rm r}}^{\rm r}$ to approximate the manifold of Riesz representers
\begin{equation}
\label{eq:manifold_test_space}
\mathfrak{M}_{\rm test}^{\rm r} =\{  \boldsymbol{\psi}_{\mu}^{\rm r} : \mu \in \mathcal{P} \},
\quad {\rm where} \;\;
(\boldsymbol{\psi}_{\mu}^{\rm r}, \mathbf{v}) = 
\mathcal{R}_{\mu}^{\rm eq,r} (  \widehat{\mathbf{u}}_{\mu},  \mathbf{v} , {\Phi}_{\mu}  ),
\;\; \forall \; \mathbf{v} \in \mathfrak{U}_0.
\end{equation}
Given the training set of parameters $\mathcal{P}_{\rm train,r} :=\{  \mu_{\rm r}^k \}_{k=1}^{n_{\rm train,r}} \subset \mathcal{P}$, we compute the solution to the ROM $\widehat{\mathbf{u}}_{\mu}$ for all $\mu \in \mathcal{P}_{\rm train,r}$ and then the Riesz representation of the residual $\boldsymbol{\psi}_{\mu}^{\rm r}$ using \eqref{eq:manifold_test_space}; 
finally, we use POD to construct the test space $\mathfrak{Y}_{J_{\rm r}}^{\rm r}$. 

On the other hand, we choose the weights $\boldsymbol{\rho}^{\rm eq,r} $ so that (i) $\boldsymbol{\rho}^{\rm eq,r} $ is as sparse as possible, (ii) all weights are non-negative, and (iii)
\begin{equation}
\label{eq:accuracy_constraint_res}
\| \widehat{\mathbf{R}}_{\mu}^{\rm eq,r} - \widehat{\mathbf{R}}_{\mu}^{\rm hf,r}   \|_2 \ll 1,
\quad
{\rm with} \;\;
\widehat{\mathbf{R}}_{\mu}^{\rm hf,r} = \sum_{k=1}^{N_{\rm e}}
\,  \widehat{\mathbf{R}}_{k,\mu}^{\rm r},
\end{equation} 
for all $\mu \in \mathcal{P}_{\rm train,r}$. Proceeding as in section \ref{sec:EQP}, we find that
$\| \widehat{\mathbf{R}}_{\mu}^{\rm eq,r} - \widehat{\mathbf{R}}_{\mu}^{\rm hf,r}   \|_2 = 
\| \mathbf{G}_{\mu}^{\rm r} \,  \boldsymbol{\rho}^{\rm eq,r}  - \widehat{\mathbf{R}}_{\mu}^{\rm hf,r}   \|_2$ with 
$\mathbf{G}_{\mu}^{\rm r} = [\widehat{\mathbf{R}}_{1,\mu}^{\rm r},\ldots,\widehat{\mathbf{R}}_{N_{\rm e},\mu}^{\rm r} ]$:
\begin{equation}
\label{eq:sparse_representation_res}
\min_{  \boldsymbol{\rho} \in \mathbb{R}^{N_{\rm e}} }
\;
\| \boldsymbol{\rho}   \|_{\ell^0},
\quad
{\rm s.t} \quad
\left\{
\begin{array}{l}
\|\mathbf{G}^{\rm r} \boldsymbol{\rho} - \mathbf{b}^{\rm r}   \|_{\star} \leq \delta; \\[3mm]
\boldsymbol{\rho} \geq \mathbf{0}; \\
\end{array}
\right.
\; \;
\mathbf{G}^{\rm r}   = \left[
\begin{array}{l}
\mathbf{G}_{\mu_{\rm r}^1}^{\rm r}   \\
\vdots \\
\mathbf{G}_{\mu_{\rm r}^{n_{\rm train,r}}}^{\rm r}   \\
\end{array}
\right]
\;\;
\mathbf{b}^{\rm r}  = \left[
\begin{array}{l}
\widehat{\mathbf{R}}_{\mu_{\rm r}^1}^{\rm hf,r} \\
\vdots  \\
\widehat{\mathbf{R}}_{\mu_{\rm r}^{n_{\rm train,r} } }^{\rm hf,r} \\
\end{array}
\right].
\end{equation}
Similarly to \eqref{eq:sparse_representation}, approximate solutions to \eqref{eq:sparse_representation_res} can be obtained using the non-negative least-squares method.

Algorithm \ref{alg:off_on_residual} summarizes the offline  procedure to compute the empirical test space $\mathfrak{Y}_{J_{\rm r}}^{\rm r}$ and the weights 
$\boldsymbol{\rho}^{\rm eq,r}$, and the online procedure to rapidly compute $\widehat{\mathfrak{R}}_{\mu}^{\rm eq} $. Note that the offline procedure takes as input the ROM and the tolerances $tol_{\rm es}$ (for $\mathfrak{Y}_{J_{\rm r}}^{\rm r}$) and  $tol_{\rm eq,r}$ (for $\boldsymbol{\rho}^{\rm eq,r}$).
In \ref{sec:appendix_remarks} (cf. Remark \ref{remark:justification_residual}), we present a rigorous justification of our approach.

\begin{algorithm}[H]                      
\caption{Offline-online estimation of the dual residual norm}     
\label{alg:off_on_residual}     

\textbf{Offline stage}
\medskip

\begin{algorithmic}[1]
\State
Given $\mathcal{P}_{\rm train,r} = \{  \mu_{\rm r}^k \}_{k=1}^{n_{\rm train,r}} $, estimate $\widehat{\mathbf{u}}_{\mu}$ by solving \eqref{eq:GalROM} for all $\mu \in \mathcal{P}_{\rm train,r}$.
\medskip

\State
Compute the Riesz representer
$\boldsymbol{\psi}_{\mu}^{\rm r}$ using \eqref{eq:manifold_test_space}
for all $\mu \in \mathcal{P}_{\rm train,r}$.
\medskip

\State
$[\boldsymbol{\eta}_1,\ldots,\boldsymbol{\eta}_{J_{\rm r}}] = $ \texttt{POD} $( 
\{ \boldsymbol{\psi}_{\mu_{\rm r}^k  }^{\rm r} \}_k, (\cdot,\cdot), tol_{\rm es}    )$.
\medskip

\State
Compute  $\mathbf{G}^{\rm r},\mathbf{b}^{\rm r}$ in \eqref{eq:sparse_representation_res}.
\medskip

\State
$[\boldsymbol{\rho}^{\rm eq,r}] =  \texttt{lsqnonneg} \left( \mathbf{G}^{\rm r}, \mathbf{b}^{\rm r}, tol_{\rm eq,r} \right).$
\medskip
\end{algorithmic}
\medskip

\textbf{Online stage} (for any given $\mu \in \mathcal{P}$ and estimate 
$\widehat{\boldsymbol{\alpha}}_{\mu}$)
\medskip 

\begin{algorithmic}[1]
\State
Compute the local residuals $\{  \widehat{\mathbf{R}}_{k,\mu}^{\rm r}  : k \in \texttt{I}_{\rm eq,r} \}$
using \eqref{eq:residual_calculation_b}.
\medskip

\State
Compute 
$ \widehat{\mathbf{R}}_{\mu}^{\rm eq,r} $ using 
\eqref{eq:residual_calculation_a} and return 
$
\widehat{\mathfrak{R}}_{\mu}^{\rm eq} = \|  \widehat{\mathbf{R}}_{\mu}^{\rm eq,r}  \|_2$.
\end{algorithmic}
\end{algorithm}

\subsection{Model problem: potential flow past a parameterized airfoil}
\label{sec:model_problem}

We introduce the domain 
$\Omega_{\mu} = 
\Omega_{\rm box} \setminus \Omega_{\rm naca,\mu}
 \subset \mathbb{R}^2$ such that  $\Omega_{\rm box} = (x_{\rm min},x_{\rm max})\times (-H,H)$,
 $x_{\rm min} = -2, x_{\rm max} = 6, H=4$, and 
\begin{equation}
\label{eq:naca_airfoil}
\Omega_{\rm naca,\mu} = \left\{
{x} \in (0,1)^2 : 
- f_{\rm naca}(x_1, \mu_1) 
< x_2
< f_{\rm naca}(x_1, \mu_2) 
\right\} ,
\end{equation}
where $f_{\rm naca}(s, \texttt{th}) = 5 \texttt{th} \left( 
0.2969 \sqrt{s} -0.1260 s - 0.3516 s^2 + 0.2843 s^3 -0.1036 s^4 \right)$, and 
$\mu_1,\mu_2$ are positive parameters. 
Then, we introduce  problem \eqref{eq:model_problem} with Dirichlet datum  $h$ given by
\begin{subequations}
\label{eq:model_problem_specialized}
\begin{equation}
h_{\mu}({x})
\; = \;
\left\{
\begin{array}{ll}
0 &  {\rm on} \; \Gamma_{\rm btm} = (x_{\rm min}, x_{\rm max}) \times\{ -H \} \\[3mm]
1  &  {\rm on} \; \Gamma_{\rm top} = (x_{\rm min}, x_{\rm max}) \times\{ H \}, \\[3mm]
\displaystyle{ \frac{x_2+H}{2H} }  &  {\rm on} \; \Gamma_{\rm out} = \{ x_{\rm max} \}  \times (-H, H)  , \\[3mm]
\displaystyle{ 
\frac{\bar{h}_{\mu}( \frac{x_2+H}{2H}    ) -  \bar{h}_{\mu}(0)  }{
\bar{h}_{\mu}(1) -  \bar{h}_{\mu}(0)}}  & 
 {\rm on} \; \Gamma_{\rm in} = 
\{ x_{\rm min} \}  \times (-H, H) ,       \\[3mm]
1   & 
 {\rm on} \;    \partial \Omega_{\rm naca,\mu} , \\
\end{array}
\right.
\end{equation}
with
\begin{equation}
 \bar{h}_{\mu}(t) 
 =\frac{1}{2} \left(
1 +\frac{1}{\pi} 
\arctan \left( 10 (t - \mu_3)  \right)
 +\frac{1}{\pi} 
\arctan \left( 10 (t - \mu_4)  \right)
\right).
\end{equation}
Finally, we shall define the vector of parameters
$\mu = [\mu_1,\mu_2,\mu_3,\mu_4]$: in the following we assume that $\mu$ belongs to the compact region
\begin{equation}
\mathcal{P} =  [0.09,0.15]^2 \times [0.1,0.3]\times [0.6,0.8].
\end{equation}
\end{subequations}

\begin{figure}[h!]
\centering
 \subfloat[ ] 
{  \includegraphics[width=0.48\textwidth]
 {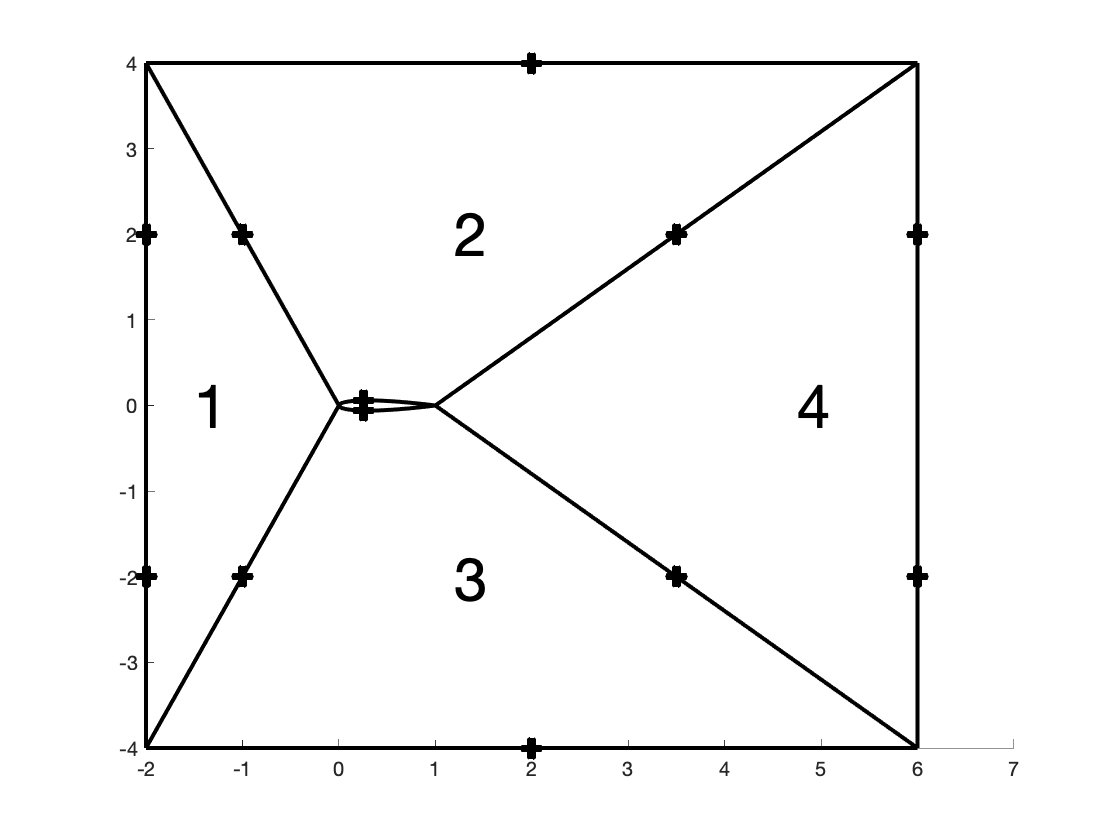}}
    ~~
\subfloat[ ] 
{  \includegraphics[width=0.48\textwidth]
 {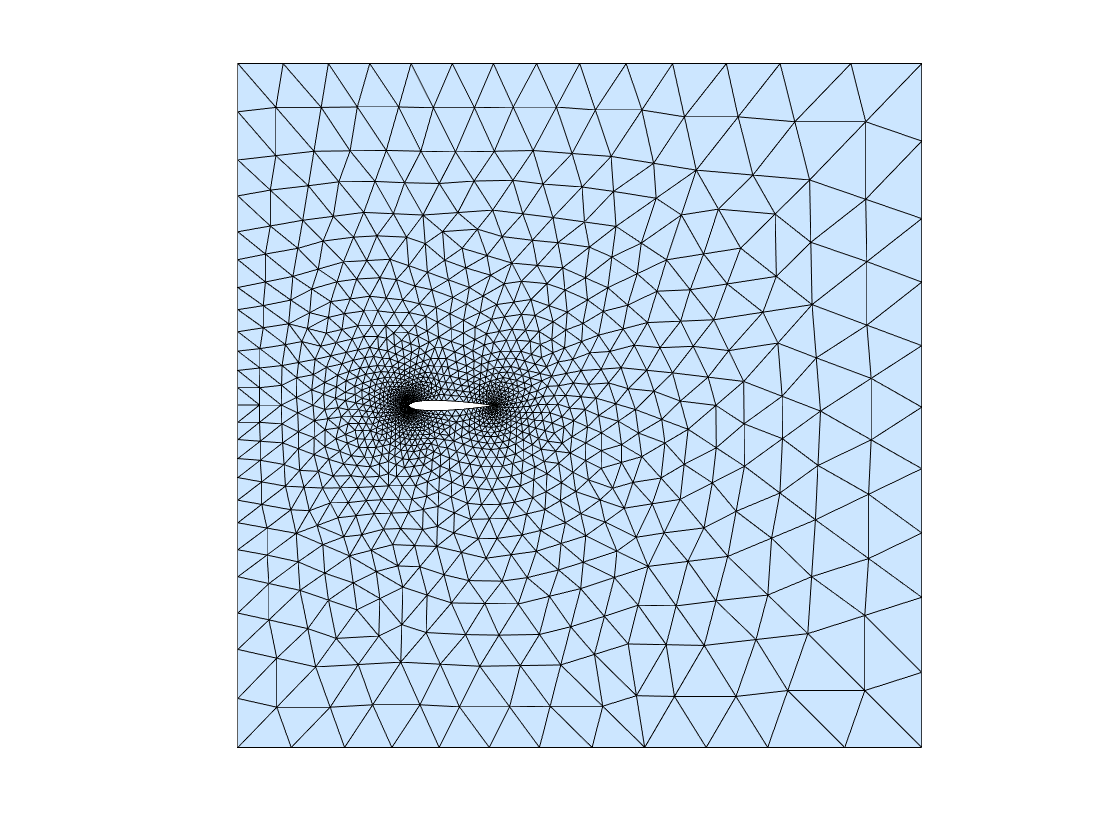}}
  
\caption{Model problem.
(a) coarse-grained partition associated with ${\Phi}$,
(b) computational mesh.}
\label{fig:vis_mesh}
\end{figure}

In order to define the mapping ${\Phi}$, 
(i)
we partition $\Omega_{\mu}$ into the four regions $\{ \Omega_{i,\mu}  \}_{i=1}^4$ depicted in Figure \ref{fig:vis_mesh}(a) and we define
$\Omega := \Omega_{\bar{\mu}}$ and
$\{ \Omega_{i} = \Omega_{i,\bar{\mu}}  \}_{i=1}^4$ with  $\bar{\mu} = [0.12,0.12,0,0]$;
(ii) we define the parameterized  Gordon-Hall maps (see \cite{gordon1973construction}) 
${\Psi}_{{\rm gh}, i}: (0,1)^2 \times \mathcal{P} \to \Omega_{i,\mu}$ for $i=1,\ldots,4$; and
(iii) we define ${\Phi}$ such that
\begin{equation}
\label{eq:Phi}
{\Phi}_{\mu}({x}) \, = \,
{\Psi}_{{\rm gh}, i} \left(
{\Psi}_{{\rm gh}, i}^{-1}
\left( {x}, \bar{\mu} \right),  \mu \right),
\quad
{\rm if} \; 
{x} \in \Omega_i,
\;\;i=1,2,3,4;
\end{equation}
note that 
${\Phi}_{\mu}$ coincides with  the identity map 
$\texttt{id}$
for $\mu=\bar{\mu}$ and 
${\Phi}_{\mu} \equiv  \texttt{id}$ in $\Omega_1 \cup \Omega_4$ for all $\mu \in \mathcal{P}$.
We remark that several other approaches might be considered to construct exact or approximate mappings for parameterized geometries: we refer to the pMOR literature for a thorough overview (see in particular \cite{lassila2014model,rozza2021basic}).

Figure \ref{fig:vis_mesh}(b) shows the mesh 
used in our implementation. Note that the mesh is not conforming  with the coarse-grained partition. To efficiently evaluate \eqref{eq:Phi}, during the offline  stage, we compute 
$\texttt{I}_{\Phi} \in \{ 1,2,3,4 \}^{N_{\rm hf}}$ and $\{  {x}_j^{\rm hf,ref} \}_{j=1}^{N_{\rm hf}} \subset [0,1]^2$ such that
$( \texttt{I}_{\Phi}  )_j$ denotes the label of the region to which the $j$-th node of the mesh ${x}_j^{\rm hf}$ belongs, and 
${x}_j^{\rm hf,ref} :=
{\Psi}_{{\rm gh}, ( \texttt{I}_{\Phi}  )_j    }^{-1}
\left( {x}_j^{\rm hf}, \bar{\mu} \right)$.
Then, given a new value of the parameter $\mu \in \mathcal{P}$, we simply compute the mapped nodes using the identity
\begin{equation}
\label{eq:onlinePhi}
{\Phi}_{\mu}( {x}_j^{\rm hf}     )
=
{\Psi}_{{\rm gh}, ( \texttt{I}_{\Phi}  )_j    } 
\left( 
{x}_j^{\rm hf,ref} , \mu \right),
\quad
j=1,\ldots,N_{\rm hf}.
\end{equation}
 
\subsection{Numerical results}
\label{sec:numerics}
We consider a $\texttt{p}=3$ FE discretization with $N_{\rm e} = 2816$ elements ($N_{\rm hf} = 12840$). To train our model, we consider $n_{\rm train} = 200$ snapshots associated with randomly-sampled parameters $\{  \mu^k \}_{k=1}^{n_{\rm train}}$. We compute the EIM approximation using a tolerance $tol_{\rm eim}=10^{-14}$: for this choice of the tolerance, we find an expansion with $M=21$ terms. In order to construct the dual-residual estimator, we consider $n_{\rm train,r}=100$ randomly-selected parameters, and we set $tol_{\rm es}=10^{-4}$ and $tol_{\rm eq,r}=10^{-10}$ in Algorithm \ref{alg:off_on_residual}. We assess performance based on $n_{\rm test}=20$ out-of-sample parameters 
$\{ \mu_{\rm t}^i  \}_{i=1}^{n_{\rm test}}$.

Figure \ref{fig:performance}(a) shows the average  out-of-sample relative error
\begin{equation}
\label{eq:Eavg}
E_{\rm avg} := \frac{1}{n_{\rm test}} 
\sum_{i=1}^{n_{\rm test}}
\;
\frac{\|  \mathbf{u}_{ \mu_{\rm t}^i }^{\rm hf} 
- \widehat{\mathbf{u}}_{ \mu_{\rm t}^i } 
 \|    }{ \|  \mathbf{u}_{ \mu_{\rm t}^i }^{\rm hf}   \| 
 }
\end{equation}
for various choices of the size $N$ of the reduced{-order} basis  $\mathbf{Z}$ and for two different tolerances $tol_{\rm eq}$ for $\boldsymbol{\rho}^{\rm eq}$.
We observe that the ROM reaches a plateau that depends on the error $\max_{ \mu \in \mathcal{P}}  \| \mathcal{H} ( \mathbf{h}_{\mu} - \widehat{\mathbf{h}}_{\mu} ) \|$ associated with EIM.
Figure \ref{fig:performance}(b) shows the { percentage of sampled elements $Q/N_{\rm e} \cdot 100$}  for several choices  of $N$ and for the two tolerances considered: note that $Q$ grows linearly with $N$. 
{ Speedups of the hyper-reduced ROM compared to the ROM with HF quadrature scale with $N_{\rm e}/Q$ and thus range from $\mathcal{O}(30)$ to $\mathcal{O}(100)$ for $N=1,\ldots,6$ and $tol_{\rm eq}=10^{-10}$: 
this can be explained by observing that, for the range of $N$ and $Q$ considered and for serial implementations, the online cost is dominated by the cost of assembling the reduced system at each Newton iteration.}

\begin{figure}[H]
\centering
\subfloat[ ] {
\begin{tikzpicture}[scale=0.7]
\begin{semilogyaxis}[
xlabel = {\LARGE {$N$}},
  ylabel = {\LARGE {$E_{\rm avg}$}},
  line width=1.2pt,
  mark size=3.0pt,
  ymin=0.00001,   ymax=0.05,
  ]
    \addplot[line width=1.pt,color=black,mark=o] table {data/laplace/gal/proj.dat};\label{ROMperformance:proj_error}
     \addplot[line width=1.pt,color=magenta,mark=square]  table {data/laplace/gal/hf_quad.dat};\label{ROMperformance:hf_quad}
  \addplot[line width=1.pt,color=red,mark=diamond]  table {data/laplace/gal/tolm6.dat};\label{ROMperformance:tol1em6}
  \addplot[line width=1.pt,color=blue,mark=triangle*] table {data/laplace/gal/tolm10.dat};\label{ROMperformance:tol1em10}
  
\end{semilogyaxis}

\end{tikzpicture}
}
~~~
\subfloat[ ] {
\begin{tikzpicture}[scale=0.7]
\begin{axis}[
xlabel = {\LARGE {$N$}},
  ylabel = {\LARGE {$\%$ sampled elements}},
  line width=1.2pt,
  mark size=3.0pt,
 ymin=0,   ymax=5,
  ]
 
  \addplot[line width=1.pt,color=red,mark=diamond]  table {data/laplace/gal/EQtolm6.dat};
  \addplot[line width=1.pt,color=blue,mark=triangle*] table {data/laplace/gal/EQtolm10.dat};  
\end{axis}

\end{tikzpicture}
}

\bgroup
\sbox0{\ref{data}}%
\pgfmathparse{\ht0/1ex}%
\xdef\refsize{\pgfmathresult ex}%
\egroup
\caption[Caption in ToC]{
Potential flow past an airfoil.
Performance  with respect to $N$ and $tol_{\rm eq}$.
(a) average out-of-sample error $E_{\rm avg}$:
projection  \tikzref{ROMperformance:proj_error},
Galerkin ROM  with HF quadrature \tikzref{ROMperformance:hf_quad},
$tol_{\rm eq}=10^{-6}$ \tikzref{ROMperformance:tol1em6},
$tol_{\rm eq}=10^{-10}$ \tikzref{ROMperformance:tol1em10}.
(b) number of sampled elements $Q$ 
for two choices of $tol_{\rm eq}$ 
(see (a)).}
 \label{fig:performance}
\end{figure}

Figure \ref{fig:residual_error} investigates the performance of the error estimation procedure.
Figures \ref{fig:residual_error}(a) and (b) show the behavior of the size ${J_{\rm r}}$ of the empirical test space $\mathfrak{Y}_{J_{\rm r}}^{\rm r}$ and of the percentage $Q_{\rm r}/N_{\rm e} \cdot 100$ of sampled elements for several ROMs associated with different choices of $N$ and two EQ tolerances $tol_{\rm eq}$. Note that
${J_{\rm r}}$ and $Q$
increase rapidly for moderate values of $N$ and then reach a plateau:  interestingly,  this is in qualitative agreement with the behavior of the error shown in 
Figure \ref{fig:performance}(a).
Figure \ref{fig:residual_error}(c) shows the relationship between the (estimated) dual residual norm  and the relative error; to generate the results, we consider the ROM with $N=5$ modes and two different choices of $tol_{\rm eq}$. We observe that the relative error is approximately proportional to the dual residual norm;  furthermore, the prediction obtained using our approach is extremely accurate for all tests considered.

\begin{figure}[h!]
\centering
\subfloat[ ] {
\begin{tikzpicture}[scale=0.7]
\begin{axis}[
xlabel = {\LARGE {$N$}},
  ylabel = {\LARGE {$J_{\rm r}$} },
  line width=1.2pt,
  mark size=3.0pt,
  ymin=1,   ymax=30
  ]
  \addplot[line width=1.pt,color=red,mark=diamond]  table {data/laplace/gal/Jresm10.dat}; 
  \addplot[line width=1.pt,color=blue,mark=triangle*] table {data/laplace/gal/Jresm12.dat};
  
\end{axis}

\end{tikzpicture}
}
~~~
\subfloat[ ] {
\begin{tikzpicture}[scale=0.7]
\begin{axis}[
xlabel = {\LARGE  {$N$}},
  ylabel = {\LARGE {$\%$ sampled elements} },
  line width=1.2pt,
  mark size=3.0pt,
 ymin=0,   ymax=10,
  legend style={at={(0,0)},anchor=south west}
  ]
 
  \addplot[line width=1.pt,color=red,mark=square]  table {data/laplace/gal/Qresm10.dat};
  \addplot[line width=1.pt,color=blue,mark=triangle*] table {data/laplace/gal/Qresm12.dat};  
\end{axis}
\end{tikzpicture}
}

\subfloat[ ] {
\begin{tikzpicture}[scale=0.7]
\begin{loglogaxis}[
xlabel = {\LARGE {dual residual}},
  ylabel = {\LARGE {$H^1$ rel error}},
legend entries = {hf res vs error, estimated res vs error},
  line width=1.2pt,
  mark size=3.0pt,
 ymin=0.00001,   ymax=0.01,
legend style={at={(0.03,0.9)},anchor=west,font=\Large}
  ]
 
  \addplot[only marks,color=red,mark=square]  table {data/laplace/gal/EQres.dat};
  
  \addplot[only marks,color=blue,mark=triangle*] table {data/laplace/gal/hfres.dat};  
  
\end{loglogaxis}
\end{tikzpicture}
}

\caption{
Potential flow past an airfoil.
Dual residual norm estimation.
(a)  size of the empirical test space  ${J_{\rm r}}$ with respect to  $N$, for $tol_{\rm eq}=10^{-10}$ 
\tikzref{ROMperformance:tol1em6} 
and
$tol_{\rm eq}=10^{-12}$ 
\tikzref{ROMperformance:tol1em10}.
(b) percentage of sampled elements $Q_{\rm r}$     with respect to  $N$, for $tol_{\rm eq}=10^{-10}$ 
$tol_{\rm eq}=10^{-12}$ \tikzref{ROMperformance:tol1em10}.
(c) exact  and estimated residual vs relative $H^1$ error.
($tol_{\rm eq,r}=10^{-10}, tol_{\rm es}=10^{-4}$)}
\label{fig:residual_error}
\end{figure}  

\section{Extension to nonlinear problems}
\label{sec:nonlinear_pb}
We discuss the application of the DtM framework to vector-valued nonlinear steady PDEs. 
Towards this end, 
we denote by $D_{\rm eq}>1$ the number of equations/state variables and we introduce the FE space
${\mathcal{X}}_{\mathcal{T}}^{\rm v} = \prod_{i=1}^{D_{\rm eq}} \, {\mathcal{X}}_{\mathcal{T}}$.
Given ${\mathbf{u}} = [{\mathbf{u}}_1^T, \ldots,{\mathbf{u}}_{D_{\rm eq}}^T ]^T \in \mathbb{R}^{{N}_{\rm hf,v}}$, ${N}_{\rm hf,v} = N_{\rm hf} \cdot D_{\rm eq}$, we can define 
(i) the
corresponding FE fields ${u}_{\Phi}$ and ${u}$ by applying \eqref{eq:vector2field} to ${\mathbf{u}}_1, \ldots,{\mathbf{u}}_{D_{\rm eq}}$;
(ii)  
the unassembled field ${\mathbf{u}}^{\rm un}  \in \mathbb{R}^{n_{\rm lp}, N_{\rm e}, D_{\rm eq}}$ such that
$\left(  {\mathbf{u}}^{\rm un} \right)_{i,k,d} = \left( \mathbf{E}_k \mathbf{u}_d  \right)_{i}$;
(iii) the set of indices ${\texttt{I}}_{\rm dir}\subset \{1,\ldots,{N}_{\rm hf,v}\}$ associated with degrees of freedom specified by Dirichlet boundary conditions; 
(iv) the spaces
${\mathfrak{U}} =  (\mathbb{R}^{{N}_{\rm hf,v}}, \langle \cdot, \cdot \rangle    )$ and
${\mathfrak{U}}_0 = 
\left\{
{\mathbf{w}} \in {\mathfrak{U}}: \,
{\mathbf{w}} (  {\texttt{I}}_{\rm dir}  ) = \mathbf{0}
\right\}$, where the inner product $ \langle \cdot, \cdot \rangle $ 
and the induced norm 
$ \vertiii{\cdot} = \sqrt{\langle \cdot, \cdot \rangle }$
will be introduced in section
\ref{sec:RANS_nonlinear} (cf. \eqref{eq:inner_product_rans}); and
(v) the vector-valued restriction operator
${\mathbf{E}}_k^{\rm v}: \mathbb{R}^{{N}_{\rm hf,v}} \to \mathbb{R}^{n_{\rm lp}\cdot   D_{\rm eq} } $ such that
$$
{\mathbf{E}}_k^{\rm v}
{\mathbf{w}}
\, = \,
\left[
\begin{array}{l}
{\mathbf{E}}_k \mathbf{w}_1\\
\vdots \\
{\mathbf{E}}_k \mathbf{w}_{D_{\rm eq}} \\
\end{array}
\right],
\quad
 \, 
{\mathbf{w}}
\, = \,
\left[
\begin{array}{l}
\mathbf{w}_1\\
\vdots \\
\mathbf{w}_{D_{\rm eq}} \\
\end{array}
\right],
\;\;
k=1,\ldots,N_{\rm e}.
$$
We further introduce  the   vector-valued extension operator
${\mathcal{H}}:   {\mathfrak{H}} \to {\mathfrak{U}}$, which is a generalization of the extension operator \eqref{eq:extension_operator}: we provide further comments about its definition in
 section  \ref{sec:RANS_nonlinear}.

We consider problems of the form: given $\mu \in \mathcal{P}$, find ${\mathbf{u}}_{\mu}^{\rm hf} \in {\mathfrak{U}}$ such that
\begin{equation}
\label{eq:abstract_nonlinear_formulation}
{\mathbf{u}}_{\mu}^{\rm hf}({\texttt{I}}_{\rm dir}   )
=
{\mathbf{h}}_{\mu},
\quad
\mathcal{R}_{\mu}^{\rm hf} \left(
{\mathbf{u}}_{\mu}^{\rm hf}, \, 
{\mathbf{v}}
\right)
=
0
\quad
\forall \, {\mathbf{v}} \in {\mathfrak{U}}_0,
\end{equation}
where $\mathcal{R}^{\rm hf} : 
{\mathfrak{U}} \times 
{\mathfrak{U}}_0 \times 
\mathcal{P} \to \mathbb{R}  $ is nonlinear in the first (and possibly the third) argument and is associated with the  system of $D_{\rm eq}>1$ PDEs of interest.
Note that, for continuous discretizations of second-order elliptic PDEs,  we might write the residual $\mathcal{R}^{\rm hf}$ as the sum of local residuals associated with integration over the elements of the HF discretization:
\begin{equation}
\label{eq:elementwise_residual}
\mathcal{R}_{\mu}^{\rm hf} \left(
{\mathbf{w}}, \, 
{\mathbf{v}}
\right)
\, = \, \sum_{k=1}^{N_{\rm e}} \,
r_{\mu,k} \left(
{\mathbf{E}}_k^{\rm v} {\mathbf{w}}, \;
{\mathbf{E}}_k^{\rm v} {\mathbf{v}},
{\Phi}_{\mu}( {\texttt{X}}_k^{\rm hf}   )
\right).
\end{equation}
As in section \ref{sec:linear_case}, we denote by
${u}_{\mu}^{\rm hf} \in {\mathfrak{X}}_{\mathcal{T}}$ (and
${u}_{\mu,\Phi_{\mu}}^{\rm hf} \in {\mathfrak{X}}_{\Phi_{\mu}(\mathcal{T})}$) the HF FE solution  field associated with the parameter $\mu\in \mathcal{P}$.

In section \ref{sec:hf_solver_nonlinear}, we present the HF solver: we introduce relevant notation and we highlight the building blocks of the full-order-model implementation.
In section \ref{sec:ROM_nonlinear}, we introduce the projection-based LSPG ROM. 
In section \ref{sec:RANS_nonlinear} we introduce the model problem and the geometric parameterization; finally, in section \ref{sec:RANS_numerics}, we present results of the numerical experiments. 

\subsection{Implementation of the high-fidelity solver}
\label{sec:hf_solver_nonlinear}

We introduce the algebraic counterpart of $\mathcal{R}^{\rm hf}$,
\begin{subequations}
 \label{eq:algebraic_nonlinear}
\begin{equation}
\label{eq:algebraic_nonlinear_residual}
\boldsymbol{\mathcal{R}}^{\rm hf}: \mathbb{R}^{ {N}_{\rm hf,v}} \times \mathcal{P} \to 
\mathbb{R}^{ {N}_{\rm hf,v}},
\qquad
\left(
\boldsymbol{\mathcal{R}}_{\mu}^{\rm hf} 
\left(
{\mathbf{w}} 
\right)
\right)_i
=
\mathcal{R}_{\mu}^{\rm hf} \left(
{\mathbf{w}} , \, 
{e}_i
\right),
\qquad
i=1,\ldots, {N}_{\rm hf,v},
\end{equation}
where ${e}_1,\ldots, {e}_{{N}_{\rm hf,v}}$ are the vectors of the canonical basis in $\mathbb{R}^{{N}_{\rm hf,v} }$, and its Jacobian 
$\boldsymbol{\mathcal{J}}^{\rm hf}: \mathbb{R}^{ {N}_{\rm hf,v}} \times  \mathcal{P} \to 
\mathbb{R}^{ {N}_{\rm hf,v},  {N}_{\rm hf,v}  }$ such that 
\begin{equation}
\label{eq:algebraic_jacobian}
\left(\boldsymbol{\mathcal{J}}_{\mu}^{\rm hf}(  {\mathbf{w}}   ) \right)_{i,j}
=\lim_{\epsilon \to 0}
\frac{1}{\epsilon} \left(
\mathcal{R}_{\mu}^{\rm hf} \left(
{\mathbf{w}}+\epsilon {e}_j, \, 
{e}_i
\right)
-
\mathcal{R}_{\mu}^{\rm hf} \left(
{\mathbf{w}} , \, 
{e}_i
\right)
\right),
\; 
i,j=1,\ldots, {N}_{\rm hf,v}.
\end{equation}
Exploiting the previous definitions, we can rewrite \eqref{eq:abstract_nonlinear_formulation} as 
\begin{equation}
\label{eq:algebraic_nonlinear_pb}
{\mathbf{u}}_{\mu}^{\rm hf}({\texttt{I}}_{\rm dir}   )
=
{\mathbf{h}}_{\mu},
\quad
\boldsymbol{\mathcal{R}}_{\mu}^{\rm hf} 
\left( {\mathbf{u}}_{\mu}^{\rm hf} \right)
\,
\left( {\texttt{I} }_{\rm in}  \right)
=
\mathbf{0},
\end{equation}
with 
$ {\texttt{I} }_{\rm in} = \{1,\ldots,{N}_{\rm hf,v} \} \setminus {\texttt{I} }_{\rm dir}$.
\end{subequations}

Problem \eqref{eq:algebraic_nonlinear_pb} can be tackled using variants of the Newton's method: 
for $k=1,\ldots$ until convergence,
\begin{equation}
\label{eq:newton_nonlinear}
{\mathbf{u}}_{\mu}^{\rm hf, (k+1)}({\texttt{I}}_{\rm in}   )
=
{\mathbf{u}}_{\mu}^{\rm hf, (k)}({\texttt{I}}_{\rm in}   )
- \; 
\left( \boldsymbol{\mathcal{J}}_{\mu}^{\rm hf,(k)} ( {\texttt{I}}_{\rm in} , {\texttt{I}}_{\rm in}     ) \right)^{-1} 
\, \left(
\boldsymbol{\mathcal{R}}_{\mu}^{\rm hf,(k)}( {\texttt{I}}_{\rm in}  )
\right),
\end{equation}
with 
$\boldsymbol{\mathcal{J}}_{\mu}^{\rm hf,(k)}= \boldsymbol{\mathcal{J}}_{\mu}^{\rm hf}
\left(   {\mathbf{u}}_{\mu}^{\rm hf, (k) }     \right)$ and 
$\boldsymbol{\mathcal{R}}_{\mu}^{\rm hf,(k)}= \boldsymbol{\mathcal{R}}_{\mu}^{\rm hf}
\left(   {\mathbf{u}}_{\mu}^{\rm hf, (k)  }    \right)$. This implies that a typical HF solver for \eqref{eq:abstract_nonlinear_formulation} involves (i) an efficient sparse assembly routine to compute 
$\boldsymbol{\mathcal{R}}_{\mu}^{\rm hf,(k)}, \boldsymbol{\mathcal{J}}_{\mu}^{\rm hf,(k)}$ at each iteration, and
(ii) a direct or iterative sparse linear solver to solve \eqref{eq:newton_nonlinear} at each iteration $k$.

Exploiting \eqref{eq:elementwise_residual}, 
given $\mu \in \mathcal{P}$ and 
${\mathbf{w}}  \in {\mathfrak{U}}$, 
the assembly routine relies on local integration to build the tensors
$\mathbf{R}_{\mu}^{\rm un}, \mathbf{J}_{\mu}^{\rm un}$ such that 
\begin{subequations}
\label{eq:tensors_nonlinear}
\begin{equation}
\left( \mathbf{R}_{\mu}^{\rm un} ( {\mathbf{w}} )   \right)_{i,k,d}
=
r_{\mu,k} \left(
{\mathbf{E}}_k^{\rm v}  {\mathbf{w}}, 
{e}_{i+(d-1) n_{\rm lp} }, 
{\Phi}_{\mu}( {\texttt{X}}_k^{\rm hf} )
\right),
\quad
\left\{
\begin{array}{l}
i=1,\ldots,n_{\rm lp} \\
d=1,\ldots,D_{\rm eq} \\
k=1,\ldots,N_{\rm e} \\
\end{array}
\right.
\end{equation}
and
\begin{equation}
\left( \mathbf{J}_{\mu}^{\rm un} ( {\mathbf{w}} )   \right)_{i,i',k,d,d'}
=
\frac{\partial}{\partial  {\mathbf{w}}_{i',k,d'}^{\rm un} } 
r_{\mu,k} \left(
{\mathbf{E}}_k^{\rm v}  {\mathbf{w}}, 
{e}_{i+(d-1) n_{\rm lp} }, 
{\Phi}_{\mu}( {\texttt{X}}_k^{\rm hf} )
\right),
\quad
\left\{
\begin{array}{l}
i,i'=1,\ldots,n_{\rm lp} \\
d,d'=1,\ldots,D_{\rm eq} \\
k=1,\ldots,N_{\rm e} \\
\end{array}
\right.
\end{equation}
\end{subequations}
and then on a sparse assembly routine to define the global vector and matrix $\boldsymbol{\mathcal{R}}_{\mu}^{\rm hf} ( {\mathbf{w}}   )  , \boldsymbol{\mathcal{J}}_{\mu}^{\rm hf}( {\mathbf{w}}  )$ in \eqref{eq:algebraic_nonlinear_residual} and \eqref{eq:algebraic_jacobian}.

Algorithm \ref{alg:nonlinear_assembly} outlines the procedure implemented in a typical FE assembler. The function 
$\texttt{local}{\_}\texttt{assembler}$ computes the tensors defined in  \eqref{eq:tensors_nonlinear}: note that the function takes as input 
(i)
the restriction of the FE field to the $k$-th element
$\left(  {\mathbf{w}}^{\rm un} \right)_{\cdot, k, \cdot}$,
(ii)
the nodes of the $k$-th element
${\Phi}_{\mu}( {\texttt{X}}_k^{\rm hf})$,  
(iii)
a structure
$\texttt{ref}$  that contains quantities associated with the reference domain --- quadrature rule $\{ ( \omega_q^{\rm g} , {X}_q^{\rm g}    ) \}_{q=1}^{n_{\rm q}}$, evaluations of the shape functions and their gradient  in the quadrature points 
$\{  {X}_q^{\rm g}    \}_{q=1}^{n_{\rm q}}$
 --- and 
(iv) 
 the structure $\texttt{input}{\_}\texttt{data}$ that contains inputs associated with the particular problem at hand --- e.g., functional expression of the source term, coefficients.
 Given  the local tensors $\mathbf{R}_{\mu}^{\rm un}$ and $\mathbf{J}_{\mu}^{\rm un}$, we define the corresponding indices associated with the global enumeration (cf. Lines 4 and 6 of the Algorithm) and then we assemble the global vector and matrix by accumulation: 
the functions 
  $\texttt{accumarray}$ and 
$\texttt{sparse}$ are the Matlab routines used to create vectors and sparse matrices by accumulation --- analogous routines can be found in any software for scientific computing.

\begin{algorithm}[H]                      
\caption{High-fidelity assembly. 
$(\mu, {\mathbf{w}} ) \rightarrow \left( \boldsymbol{\mathcal{R}}_{\mu}^{\rm hf} ( {\mathbf{w}}   )  , \boldsymbol{\mathcal{J}}_{\mu}^{\rm hf}( {\mathbf{w}}  ) \right)$.}     
\label{alg:nonlinear_assembly}     

\begin{algorithmic}[1]

\For{$k=1,\ldots,N_{\rm e}$}
\State
$\left[
\left(   
\mathbf{R}_{\mu}^{\rm un} ( {\mathbf{w}} )  
\right)_{\cdot,k,\cdot}, 
\left(   
\mathbf{J}_{\mu}^{\rm un} ( {\mathbf{w}} )  
\right)_{\cdot,\cdot,k,\cdot,\cdot} 
\right]
=
\texttt{local}{\_}\texttt{assembler}
\left(
\left(    \mathbf{w}^{\rm un} \right)_{\cdot, k, \cdot}, 
\;
{\Phi}_{\mu}( {\texttt{X}}_k^{\rm hf}), 
\texttt{ref},  \texttt{input}{\_}\texttt{data}
\right)
$
\EndFor
\smallskip

\State
Define $\texttt{ii}_{i,k,d}^{\rm r} = \texttt{T}_{i,k} + N_{\rm hf} (d-1)$ for $i=1,\ldots,n_{\rm lp}$, $k=1,\ldots,N_{\rm e}$, $d=1,\ldots,D_{\rm eq}$.
\smallskip

\State
 $\boldsymbol{\mathcal{R}}_{\mu}^{\rm hf} ( {\mathbf{w}}   )  =
\texttt{accumarray}\left(
\texttt{ii}^{\rm r}, \mathbf{R}_{\mu}^{\rm un} ( {\mathbf{w}} )  , {N}_{\rm hf,v}, 1
\right)
$.
\smallskip

\State
Define $\texttt{ii}_{i,i',k,d,d'}^{\rm j} = \texttt{T}_{i,k} + N_{\rm hf} (d-1)$ 
and
$\texttt{jj}_{i,i',k,d,d'}^{\rm j} = \texttt{T}_{i',k} + N_{\rm hf} (d'-1)$ 
for $i,i'=1,\ldots,n_{\rm lp}$, $k=1,\ldots,N_{\rm e}$, $d,d'=1,\ldots,D_{\rm eq}$.
\smallskip

\State
$\boldsymbol{\mathcal{J}}_{\mu}^{\rm hf} ( {\mathbf{w}}   )  =
\texttt{sparse}\left(
\texttt{ii}^{\rm j}, \texttt{jj}^{\rm j}, \mathbf{J}_{\mu}^{\rm un} ( {\mathbf{w}} )  ,  {N}_{\rm hf,v},  {N}_{\rm hf,v}
\right)
$.
\end{algorithmic}
\end{algorithm}

To reduce the level of intrusiveness of the pMOR procedure, and ultimately simplify the implementation, we shall reuse the local integration routine 
$\texttt{local}{\_}\texttt{assembler}$ to assemble  ROM matrices and vectors: in the next section, we illustrate how to achieve this goal. 
We  observe that for several problems in computational mechanics, the naive Newton's method in \eqref{eq:newton_nonlinear} might not converge: in section \ref{sec:RANS_nonlinear}, we illustrate a possible remedy for the RANS equations, which is used in the numerical simulations.

\begin{remark}
\label{remark:DG}
Although we here focus on continuous FE discretizations, Algorithm \ref{alg:nonlinear_assembly} --- and the subsequent pMOR procedure --- might be generalized to finite volume and discontinuous Galerkin (DG) discretizations.
More in detail, for second-order problems, DG discretizations rely on local assembling routines that take as inputs the solution in the $k$-th element and in its neighbors. As a result, computation  of the ROM residual requires to store the basis functions in the sampled elements and in their neighbors. { Algorithm \ref{alg:nonlinear_assembly} might  also be extended to Taylor-Hood discretizations: the implementation should take into account the fact that the number of degrees of freedom per element differs among the different variables (e.g., pressure and velocity).
}
\end{remark}

\subsection{Least-squares Petrov-Galerkin  reduced-order model}
\label{sec:ROM_nonlinear}

For simplicity, we assume that ${\mathbf{h}}_{\mu} = {\boldsymbol{\Xi}} \boldsymbol{\Theta}_{\mu}$ with 
$\boldsymbol{\Theta}: \mathcal{P} \to  \mathbb{R}^M$, $M=\mathcal{O}(1)$; then, given the extension operator 
${\mathcal{H}}:   {\mathfrak{H}} \to {\mathfrak{U}}$, we obtain the affine-in-parameter lifting 
${\mathbf{e}}_{\mu} =
{\mathcal{H}}  {\mathbf{h}}_{\mu} =
{\mathbf{W}} \boldsymbol{\Theta}_{\mu}
$ with  ${\mathbf{W}}  = {\mathcal{H}}  {\boldsymbol{\Xi}}$; the extension to the non-affine case can be addressed as in section \ref{sec:linear_case} using EIM.
Given the reduced{-order} bases 
${\mathbf{Z}} =  
  [{\boldsymbol{\zeta}}_1, \ldots, {\boldsymbol{\zeta}}_N  ] \in \mathbb{R}^{{N}_{\rm hf,v}, N}$ and
  ${\mathbf{Y}} =  
  [{\boldsymbol{\xi}}_1, \ldots, {\boldsymbol{\xi}}_J  ] \in \mathbb{R}^{{N}_{\rm hf,v}, J}$
    with 
  ${\mathbf{Z}}({\texttt{I}}_{\rm dir}, :     ) = \mathbf{0}$, 
  ${\mathbf{Y}}({\texttt{I}}_{\rm dir}, :     ) = \mathbf{0}$, we seek approximations of the form
\begin{subequations}
\label{eq:AMR_ROM_all}  
  \begin{equation}
\label{eq:AMR_ROM}
\widehat{{\mathbf{u}}}_{\mu} 
=
{\mathbf{Z}}  \widehat{\boldsymbol{\alpha}}_{\mu} 
+ {\mathbf{e}}_{\mu} ,
\quad
{\rm such \;  that} \;
\widehat{\boldsymbol{\alpha}}_{\mu} \in {\rm arg} \min_{ 
\boldsymbol{\alpha} \in \mathbb{R}^N }
\;
\sup_{{\boldsymbol{\xi}} \in {\rm col} ( {\mathbf{Y}}   )   }
\;
\frac{ \mathcal{R}_{\mu}^{\rm eq}( 
{\mathbf{Z}}  \boldsymbol{\alpha}
+ {\mathbf{e}}_{\mu}, \; {\boldsymbol{\xi}} 
)     }{\vertiii{ {\boldsymbol{\xi}}   }},
 \end{equation}
  where $\mathcal{R}_{\mu}^{\rm eq}$ is the empirical weighted  residual
  \begin{equation}
\label{eq:elementwise_weighted_residual}
\mathcal{R}_{\mu}^{\rm eq} \left(
{\mathbf{w}}, \, 
{\mathbf{v}}
\right)
\, = \, \sum_{k=1}^{N_{\rm e}} \, \rho_k^{\rm eq} \; 
r_{\mu,k} \left(
{\mathbf{E}}_k^{\rm v}  {\mathbf{w}}, \;
{\mathbf{E}}_k^{\rm v}  {\mathbf{v}},
{\Phi}_{\mu}( \texttt{X}_k^{\rm hf}   )
\right);
\end{equation}
  with  $\boldsymbol{\rho}^{\rm eq} \in \mathbb{R}_+^{N_{\rm e}}$,
$\| \boldsymbol{\rho}^{\rm eq}  \|_{\ell^0} = Q \ll N_{\rm e}$. Provided that ${\boldsymbol{\xi}}_1, \ldots, {\boldsymbol{\xi}}_J$ are orthonormal in ${\mathfrak{U}}$, we can rewrite \eqref{eq:AMR_ROM} as
\begin{equation}
\label{eq:AMR_ROM_algebraic}
\widehat{\boldsymbol{\alpha}}_{\mu} \in {\rm arg} \min_{ 
\boldsymbol{\alpha} \in \mathbb{R}^N }
\;
\| \widehat{\mathbf{R}}_{N,J}( \boldsymbol{\alpha}, \mu   )  \|_2,
\;\;
\left(
\widehat{\mathbf{R}}_{N,J}( \boldsymbol{\alpha}, \mu   ) 
\right)_j =
\mathcal{R}_{\mu}^{\rm eq}( 
{\mathbf{Z}}  \boldsymbol{\alpha}
+ {\mathbf{e}}_{\mu}, \; {\boldsymbol{\xi}}_j 
), 
\end{equation}
for $j=1,\ldots,J$.
\end{subequations}
{ Problem \eqref{eq:AMR_ROM_algebraic} is a nonlinear least-squares problem: we can thus resort to the Gauss-Newton algorithm to efficiently compute the solution. We rely on the dataset $\{  ( \mu^k, \boldsymbol{\alpha}_{\mu^k}^{\rm train}   )  \}_{k=1}^{n_{\rm train}}$ (cf. (\ref{eq:accuracy_constraint_nonlinear}b)) to train a regressor $\mu \in \mathcal{P}  \mapsto \widehat{\boldsymbol{\alpha}}_{\mu}^{(0)}$ that is used to initialize the iterative Gauss-Newton scheme: we refer to \cite[section 4.1]{taddei2020space} for a thorough description of  the multi-target regression algorithm employed in this work.}

Similarly to the linear case, we should discuss the choice of (i) the trial reduced{-order basis}  ${\mathbf{Z}}$, (ii) the empirical weights $\boldsymbol{\rho}^{\rm eq}$, 
(iii) the test reduced{-order basis}  ${\mathbf{Y}}$, and the construction of the approximate dual residual estimator.
Before addressing these tasks, we derive explicit expressions for 
$\widehat{\mathbf{R}}_{N,J}( \boldsymbol{\alpha}, \mu   )$ and the Jacobian
$\widehat{\mathbf{J}}_{N,J}( \boldsymbol{\alpha}, \mu   ) = \nabla_{ \boldsymbol{\alpha}  } \widehat{\mathbf{R}}_{N,J}( \boldsymbol{\alpha}, \mu   )  $.

We denote by ${\mathbf{Z}}^{\rm un} \in \mathbb{R}^{n_{\rm lp}, N_{\rm e}, D_{\rm eq}, N}$ and 
${\mathbf{Y}}^{\rm un} \in \mathbb{R}^{n_{\rm lp}, N_{\rm e}, D_{\rm eq}, J}$ the unassembled trial and test reduced spaces, and we define
$\texttt{I}_{\rm eq} = \{ k \in \{ 1,\ldots,N_{\rm e} \} : \rho_k^{\rm eq}>0   \}$.  Then, it is easy to verify that
\begin{equation}
\label{eq:assembler_ROM_nonlinear}
\left\{
\begin{array}{ll}
\displaystyle{
\left(
\widehat{\mathbf{R}}_{N,J}( \boldsymbol{\alpha}, \mu   ) 
\right)_j
=
\sum_{k\in \texttt{I}_{\rm eq}} \; \rho_k^{\rm eq} \, \left(
\sum_{i=1}^{n_{\rm lp}}
\sum_{d=1}^{D_{\rm eq}}
\left( {\mathbf{Y}}^{\rm un} \right)_{i,k,d,j}
\left( {\mathbf{R}}^{\rm un} ( \widehat{{\mathbf{u}}}_{\mu}  )   \right)_{i,k,d}
\right)
}
&
\\[3mm]
\displaystyle{
\left(
\widehat{\mathbf{J}}_{N,J}( \boldsymbol{\alpha}, \mu   ) 
\right)_{j,n}
=
\sum_{k\in \texttt{I}_{\rm eq}} \; \rho_k^{\rm eq} \, \left(
\sum_{i,i'=1}^{n_{\rm lp}}
\sum_{d,d'=1}^{D_{\rm eq}}
\left( {\mathbf{Y}}^{\rm un} \right)_{i,k,d,j}
\left( {\mathbf{Z}}^{\rm un} \right)_{i',k,d',n}
\left( {\mathbf{J}}^{\rm un} ( \widehat{{\mathbf{u}}}_{\mu}  )   \right)_{i,i',k,d,d'}
\right)
}
&
 \\ 
\end{array}
\right.
\end{equation}
for $j=1,\ldots,J$, $n=1,\ldots,N$.

Algorithm \ref{alg:nonlinear_ROM_assembly} summarizes the assembly procedure for the computation of the reduced residual. Note that Algorithm \ref{alg:nonlinear_ROM_assembly} exploits the same local integration routines employed in Algorithm \ref{alg:nonlinear_assembly}; as in the linear case, evaluation of 
$\widehat{\mathbf{R}}_{N,J}( \boldsymbol{\alpha}, \mu   )$ , $\widehat{\mathbf{J}}_{N,J}( \boldsymbol{\alpha}, \mu )$ requires the storage of the reduced{-order basis} in $Q=| \texttt{I}_{\rm eq}|$ elements of the mesh and can thus be performed in $\mathcal{O}(Q)$ operations. 
{
We highlight that Algorithm \ref{alg:nonlinear_ROM_assembly}  depends on the underlying FE model exclusively through the local assembly routine:
provided that the latter is available, implementation of the ROM assembly routine  is straightforward.
}

\begin{algorithm}[H]                      
\caption{ROM assembly 
$(\mu, \boldsymbol{\alpha} ) \rightarrow \left( 
\widehat{\mathbf{R}}_{N,J}( \boldsymbol{\alpha}, \mu   )  , \; \; \widehat{\mathbf{J}}_{N,J}( \boldsymbol{\alpha}, \mu   )  \right)$.}     
\label{alg:nonlinear_ROM_assembly}     

\begin{algorithmic}[1]
\State
{ Compute the deformed nodes
$\{ {\Phi}_{\mu}({x}_{ i,k}^{\rm hf}): i=1,\ldots,n_{\rm lp}, k \in  \texttt{I}_{\rm eq} \}$}.
\medskip

\State
Compute
$\left(  \widehat{{\mathbf{u}}}_{\mu}^{\rm un} \right)_{i, k,d}
= \sum_n \left( {\mathbf{Z}}^{\rm un} \right)_{i,k,d,n} 
\left( \boldsymbol{\alpha} \right)_{n}$,
$i=1,\ldots,n_{\rm lp}, k \in \texttt{I}_{\rm eq}, d=1,\ldots,D_{\rm eq}$.
\medskip 

\For{$k\in \texttt{I}_{\rm eq}$}
\State
$\left[
\left(   
\mathbf{R}_{\mu}^{\rm un} ( 
\widehat{{\mathbf{u}}}_{\mu}
)  
\right)_{\cdot,k,\cdot}, 
\left(   
\mathbf{J}_{\mu}^{\rm un} ( 
\widehat{{\mathbf{u}}}_{\mu}
)  
\right)_{\cdot,\cdot,k,\cdot,\cdot} 
\right]
=
\texttt{local}{\_}\texttt{assembler}
\left(
\left(  \widehat{{\mathbf{u}}}_{\mu}^{\rm un}  \right)_{\cdot, k, \cdot}, 
\;
{\Phi}_{\mu}( {\texttt{X}}_k^{\rm hf}), 
\texttt{ref},  \texttt{input}{\_}\texttt{data}  
\right)$.
\EndFor
\smallskip

\State
Compute 
$\widehat{\mathbf{R}}_{N,J}( \boldsymbol{\alpha}, \mu   )$ and $ \widehat{\mathbf{J}}_{N,J}( \boldsymbol{\alpha}, \mu   )$  using \eqref{eq:assembler_ROM_nonlinear}.
 
\end{algorithmic}
\end{algorithm}

\subsubsection{Construction of  trial and test reduced-order spaces}
Given the snapshots $\{ {\mathbf{u}}_{\mu^k}^{\rm hf}  \}_{k=1}^{n_{\rm train}}$, we apply POD to generate the reduced-order basis ${\mathbf{Z}}$:
$$
{\mathbf{Z}} =
\texttt{POD} \left(
\{ \mathring{{\mathbf{u}}}_{\mu^k}^{\rm hf}  \}_{k=1}^{n_{\rm train}}, \langle \cdot, \cdot \rangle, tol_{\rm pod}
\right),
\quad
\mathring{{\mathbf{u}}}_{\mu}^{\rm hf} :=
{{\mathbf{u}}}_{\mu}^{\rm hf} -
{{\mathbf{e}}}_{\mu} 
\quad \mu \in \mathcal{P}.
$$
On the other hand, we resort to the procedure proposed in \cite{taddei2020space} to construct the test space. For completeness, Algorithm \ref{alg:nonlinear_emp_test_space} summarizes the data compression algorithm for the construction of the test space: 
the FE vector ${\boldsymbol{\psi}}_{k,n} \in {\mathfrak{U}}_0$ on Line 2 of the Algorithm is the Riesz representer of the Fr{\'e}chet  derivative of the residual at $ {\mathbf{u}}_{\mu^k}^{\rm hf}$ applied to the $n$-th element of the reduced-order basis;
we refer to \cite{taddei2020space} for further details and for the mathematical justification of the approach for linear problems.

\begin{algorithm}[H]                      
\caption{Construction of the test space 
${\mathfrak{Y}}_J =(  {\rm col} ({\mathbf{Y}}),  \langle \cdot, \cdot \rangle ) $. }     
\label{alg:nonlinear_emp_test_space}     

\begin{algorithmic}[1]
\For{$k\in n_{\rm train}$, $n=1,\ldots,N$}
\State
Find ${\boldsymbol{\psi}}_{k,n} \in {\mathfrak{U}}_0$ such that
$
\langle 
{\boldsymbol{\psi}}_{k,n}, \, 
{\mathbf{v}}
\rangle
=
\left(
\boldsymbol{\mathcal{J}}_{\mu^k}^{\rm hf}
\left(  {\mathbf{u}}_{\mu^k}^{\rm hf}   \right)
{\boldsymbol{\zeta}}_n   \right) \cdot
{\mathbf{v}}
$
for all ${\mathbf{v}} \in {\mathfrak{U}}_0$.
\EndFor
\smallskip

\State
${\mathbf{Y}} =
\texttt{POD} \left(
\{ {\boldsymbol{\psi}}_{k,n}  \}_{k,n},  \langle \cdot, \cdot \rangle, tol_{\rm es}
\right)$.
 
\end{algorithmic}
\end{algorithm}

\subsubsection{Empirical quadrature}
We apply a variant of the empirical quadrature procedure considered in \cite{taddei2020space}. More precisely, given the snapshots $\{  {{\mathbf{u}}}_{\mu^k}^{\rm hf}  \}_{k=1}^{n_{\rm train}}$ and the reduced-order bases
${\mathbf{Z}}$ and ${\mathbf{Y}}$, we seek 
$\boldsymbol{\rho}^{\rm eq} \in \mathbb{R}_+^{N_{\rm e}}$ such that  
(i)    $\| \boldsymbol{\rho}^{\rm eq}   \|_{\ell^0}$,   is as small as possible,
(ii) the constant function is accurately approximated in $\Omega$ (see \eqref{eq:constant_function_constraint}), and
(iii) for all $\mu \in \mathcal{P}_{\rm train,eq} =\mathcal{P}_{\rm train} \cup \{ \mu^{{\rm eq,} i}  \}_{i=1}^{n_{\rm train,eq}}  \subset \mathcal{P}$,
\begin{subequations}
\label{eq:accuracy_constraint_nonlinear}
\begin{equation}
\|  
\widehat{\mathbf{R}}_{N,J}(  \boldsymbol{\alpha}_{\mu}^{\rm train}, \mu  ) 
-
\widehat{\mathbf{R}}_{N,J}^{\rm hf}(  \boldsymbol{\alpha}_{\mu}^{\rm train}, \mu  ) 
\|_2 \ll 1
\quad
\forall \; \mu \in \mathcal{P}_{\rm train,eq},
\end{equation}
where $\widehat{\mathbf{R}}_{N,J}^{\rm hf}: \mathbb{R}^N \times \mathcal{P} \to \mathbb{R}^J$ is the reduced residual associated with the choice $\rho_k^{\rm eq} = 1$ for $k=1,\ldots,N_{\rm e}$ (reduced residual with HF quadrature), while
$\boldsymbol{\alpha}_{\mu}^{\rm train}$ is chosen such that
\begin{equation}
\label{eq:accuracy_constraint_nonlinear_b}
\boldsymbol{\alpha}_{\mu}^{\rm train}
=
\left\{
\begin{array}{ll}
\displaystyle{
\left[
\langle \mathring{{\mathbf{u}}}_{\mu}^{\rm hf}, \,
{\boldsymbol{\zeta}}_1\rangle, \ldots
\langle \mathring{{\mathbf{u}}}_{\mu}^{\rm hf}, \,
{\boldsymbol{\zeta}}_N \rangle
\right]^T
} & {\rm if} \; \mu \in \mathcal{P}_{\rm train}, \\[3mm]
\displaystyle{
{\rm arg} \min_{ 
\boldsymbol{\alpha} \in \mathbb{R}^N }
\;
\| \widehat{\mathbf{R}}_{N,J}^{\rm hf}( \boldsymbol{\alpha}, \mu   )  \|_2 
} & {\rm if} \; \mu \notin \mathcal{P}_{\rm train}. \\
\end{array}
\right.
\end{equation}
Note that \eqref{eq:accuracy_constraint_nonlinear} involves the solution to $n_{\rm train,eq}$ ROMs with HF quadrature. By combining \eqref{eq:constant_function_constraint} with \eqref{eq:accuracy_constraint_nonlinear}, we obtain a sparse representation problem of the form \eqref{eq:sparse_representation}, which can be tackled using the non-negative least-squares method (cf. section \ref{sec:EQP}).
\end{subequations}

Some comments are in order. First, in \cite{taddei2020space}, we consider a slightly different accuracy constraint in which we pre-multiply $\widehat{\mathbf{R}}_{N,J}(  \boldsymbol{\alpha}_{\mu}^{\rm train}, \mu  ) 
-
\widehat{\mathbf{R}}_{N,J}^{\rm hf}(  \boldsymbol{\alpha}_{\mu}^{\rm train}, \mu  ) $ by  the transposed reduced Jacobian: in our experience, the constraint \eqref{eq:accuracy_constraint_nonlinear} performs better for moderate values of $n_{\rm train}$. Second, to improve performance of EQ, similarly to \cite[Algorithm 1]{yano2019discontinuous}, we add training points to $\mathcal{P}_{\rm train}$ by solving the ROM with HF quadrature: we empirically observe that this choice   improves performance of the hyper-reduced ROM, particularly for moderate values of $n_{\rm train}$.
{
In addition, since the solution to 
\eqref{eq:accuracy_constraint_nonlinear_b}$_2$ is significantly less expensive than an HF solve, increasing the size of $\mathcal{P}_{\rm train,eq}$ does not significantly increase offline costs.
}

\subsubsection{A posteriori error estimation}

As for the linear case, we rely on an approximate dual residual estimator:
\begin{equation}
\label{eq:dual_residual_nonlinear}
\mathfrak{R}_{\mu}^{\rm eq} ( \widehat{{\mathbf{u}}}_{\mu}    )
:=
\sup_{ {\mathbf{v}} \in {\mathfrak{Y}}_{J_{\rm r}}^{\rm r} } \;
\frac{\mathcal{R}_{\mu}^{\rm eq,r}( 
\widehat{{\mathbf{u}}}_{\mu}  ,
 {\mathbf{v}} )}{\vertiii{ {\mathbf{v}}  }   },
 \;\;
 \mathcal{R}_{\mu}^{\rm eq,r} \left(
{\mathbf{w}}, \, 
{\mathbf{v}}
\right)
\, = \, \sum_{k=1}^{N_{\rm e}} \, \rho_k^{\rm eq,r} \; 
r_{\mu,k} \left(
{\mathbf{E}}_k^{\rm v}  {\mathbf{w}}, \;
{\mathbf{E}}_k^{\rm v}  {\mathbf{v}},
{\Phi}_{\mu}( \texttt{X}_k^{\rm hf}   )
\right); 
\end{equation}
where the space 
${\mathfrak{Y}}_{J_{\rm r}}^{\rm r} \subset {\mathfrak{U}}_0$ and the weights
$\boldsymbol{\rho}^{\rm eq,r} \in \mathbb{R}_+^{N_{\rm e}}$
are computed by extending the approach presented 
in section \ref{sec:a_posteriori} to nonlinear problems --- since the extension is straightforward, we omit the details. As discussed in the pMOR literature, rigorous and rapidly-computable \emph{a posteriori} error estimators are rarely available for nonlinear PDEs: in the numerical experiments, we verify numerically that the residual estimator is highly-correlated with the reconstruction error and can thus be used to drive an adaptive sampling strategy for the optimal selection of the training parameters $\{  \mu^k \}_{k=1}^{n_{\rm train}}$.

\subsection{Model problem: RANS simulations of the flow past the Ahmed body}
\label{sec:RANS_nonlinear}

\begin{figure}[h!]
\centering

 \subfloat[] {
\begin{tikzpicture}[scale=6]
\linethickness{0.3 mm}
\draw[ultra thick]  (0.1,0.05)--(1.044,0.05)--(1.044,0.2442)--(0.8428,0.3380)--(0.1,0.3380);
\draw [ultra thick,domain=90:180,samples=300] plot ({0.1+0.1*cos(\x)}, {0.2380+0.1*sin(\x)});
\draw [ultra thick,domain=180:270,samples=300] plot ({0.1+0.1*cos(\x)}, {0.15+0.1*sin(\x)});
\draw[ultra thick]  (0,0.2380)--(0,0.1500);

\draw [blue,ultra thick,dashed,domain=335:360,samples=30] plot ({0.8428+0.2*cos(\x)}, {0.3380+0.2*sin(\x)});
\draw[blue,dashed,ultra thick]  (0.8428,0.3380)--(1.0428,0.3380);
 \coordinate [label={right:  {\Large{$\mu$}}}] (E) at (1.0428,0.3) ;
  
\draw[ultra thick]  (-0.2,0)--(1.2,0); 
 \draw[<->]  (1.1,0)--(1.1,0.05); 
  \coordinate [label={right:  {{$h_{\rm b}$}}}] (E) at (1.1,0.05) ;

  \draw[<->]  (0,0.4)--(1.044,0.4); 
  \coordinate [label={above:  {{$L$}}}] (E) at (0.522,0.4) ;

\draw[<->]  (1.2,0.05)--(1.2,0.3380); 
  \coordinate [label={right:  {{$H_{\rm c}$}}}] (E) at (1.2,0.17) ;
 
\draw[<->]  (0.8428-0.02,0.3380-0.02)--(1.044-0.02,0.2442-0.02);  
   \coordinate [label={below:  {{$\delta$}}}] (E) at (0.9234, 0.2711) ;
 
 \draw[<->]  (0.1,0.15)--(0,0.15);  
   \coordinate [label={above:  {{$R$}}}] (E) at (0.05, 0.15) ;
\end{tikzpicture}
}

\subfloat[] {
\includegraphics[width=0.48\textwidth]
{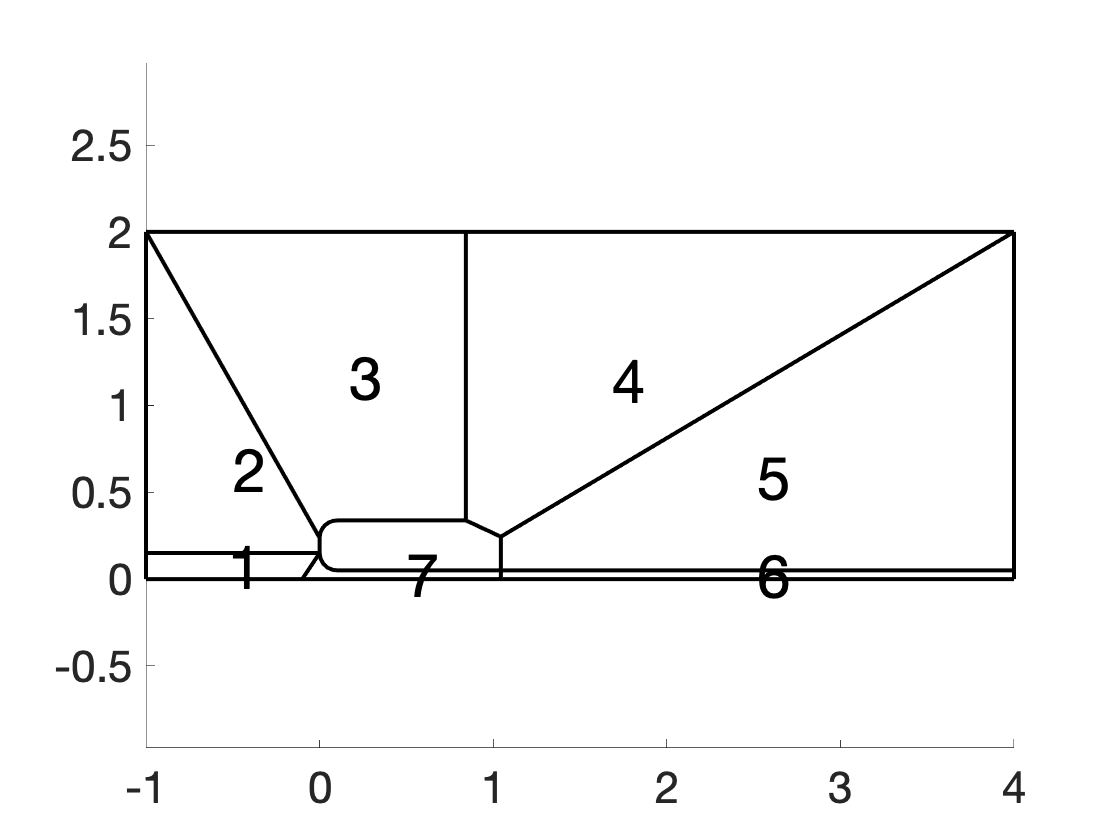}
  }
~~~~
\subfloat[] {
\includegraphics[width=0.48\textwidth]
{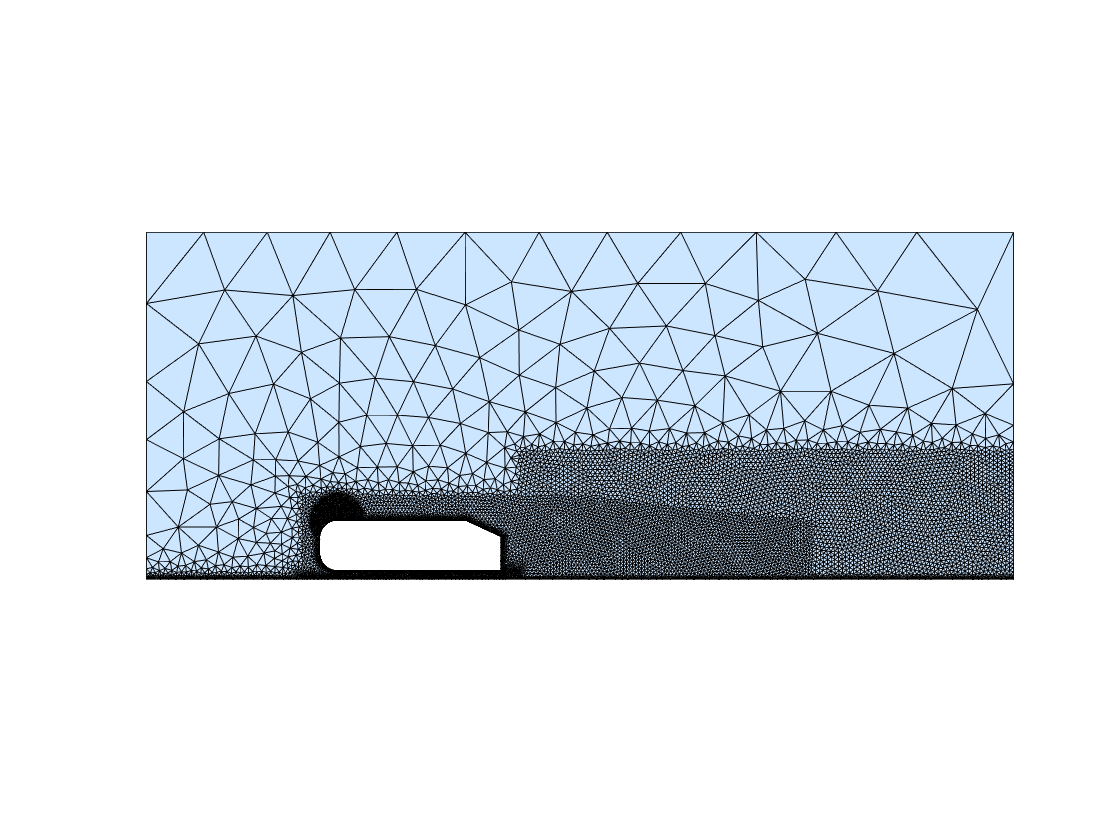}
  }
 \caption{Ahmed's body.
(a)  description of the body geometry.
(b) coarse-grained partition associated with ${\Phi}$. (c) computational mesh.
 $L=1.044$, $R=0.1$, $h_{\rm b} = 0.05$, $H_{\rm c} = 0.288$, $\delta =0.222$. 
 }
\label{fig:ahmed_body}
\end{figure}

We consider a two-dimensional flow past a parameterized Ahmed body, \cite{Ahmed1984,Lienhart2002},  at moderate Reynolds number.
The geometry of the body is prescribed in Figure \ref{fig:ahmed_body}(a).
The problem is parameterized with respect to the slant angle $\mu$: in the numerical simulations, we consider angles between $5^o$ and $50^o$ degrees, $\mathcal{P} = [5,50]$.
Reynolds number is defined as the ratio between horizontal inflow velocity times the height $H_{\rm c}$ of the body divided by kinematic viscosity; in our tests, we consider $Re=3 \cdot 10^3$. 

We consider the incompressible RANS equations with 
Spalart-Allmaras closure model 
 (see 
\cite{Spalart1992, Allmaras2012}). We denote by $U$  the set of $D_{\rm eq}=4$ state variables 
which comprises the two velocity components, pressure and turbulent viscosity. We refer to {  \ref{sec:RANS_appendix} for further details concerning the mathematical model and boundary conditions.}

We discretize the problem using  a  
SUPG FE method,  \cite{tezduyar1991stabilized,tezduyar2000finite}.
We add the Pressure-Stabilized Petrov-Galerkin (PSPG) term to eliminate  spurious modes in the pressure solution when considering the same polynomial order for pressure and velocity; furthermore, we consider the 
least-squares incompressibility constraint (LSIC)  stabilization term  to improve  accuracy and conditioning of the discrete problem,  \cite{Franca1992,Gelhard2005,braack2007stabilized}. 
In all our tests, we consider the mesh depicted in Figure \ref{fig:ahmed_body}(c) that consists of $N_{\rm e} = 26602$ P2 elements (${N}_{\rm hf,v} = 221960$).
To compute the HF solution, we resort to a pseudo-time continuation (see, e.g., \cite{kelley1998convergence}) procedure to determine an accurate initial condition for the Newton solver; then, we apply the Newton scheme 
\eqref{eq:newton_nonlinear}
with 1D  line search.

Similarly to section \ref{sec:model_problem}, we partition the domain 
$\Omega_{\mu}$ into the seven regions $\{ \Omega_{i,\mu}  \}_{i=1}^7$
depicted in Figure \ref{fig:ahmed_body}(b),  we define 
parameterized  Gordon-Hall maps in each of the seven regions, and we define ${\Phi}$ as in \eqref{eq:Phi}. Note that 
${\Phi}$ is the identity in the regions $1,2,6,7$.

\begin{figure}[h!]
\centering
\subfloat[] {
\includegraphics[width=0.48\textwidth]
{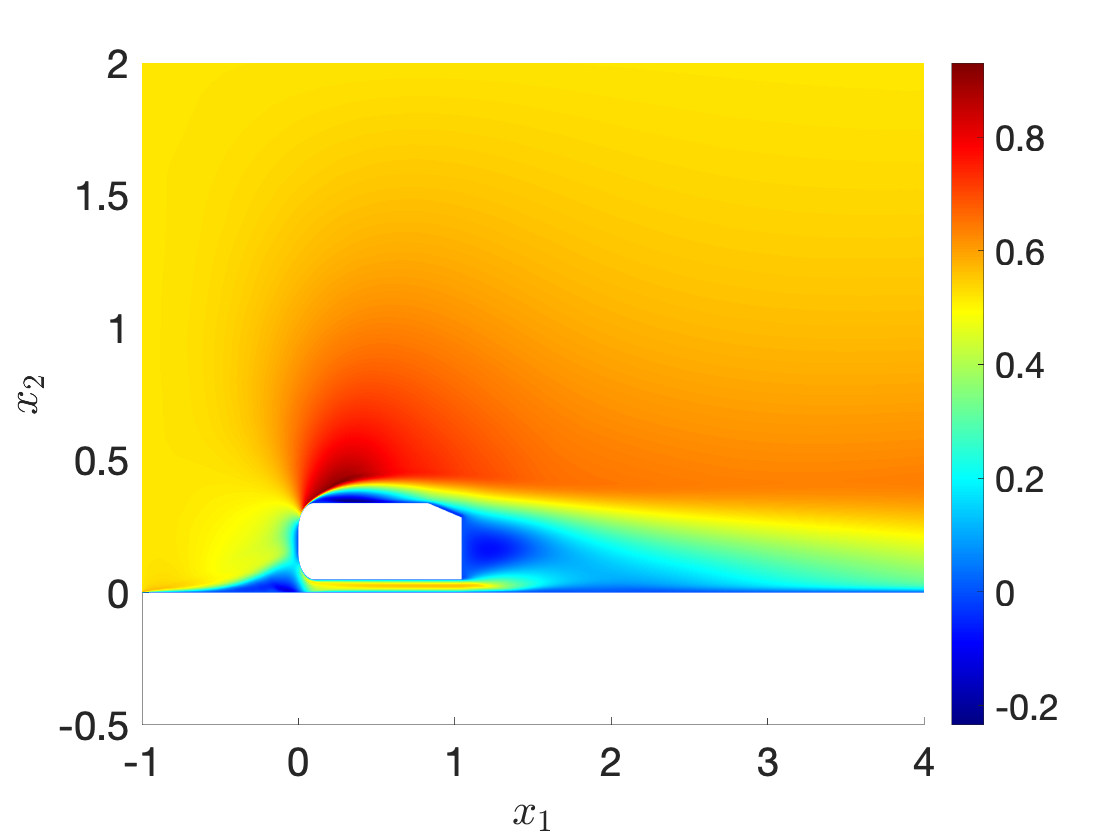}
  }
~~~~
\subfloat[] {
\includegraphics[width=0.48\textwidth]
{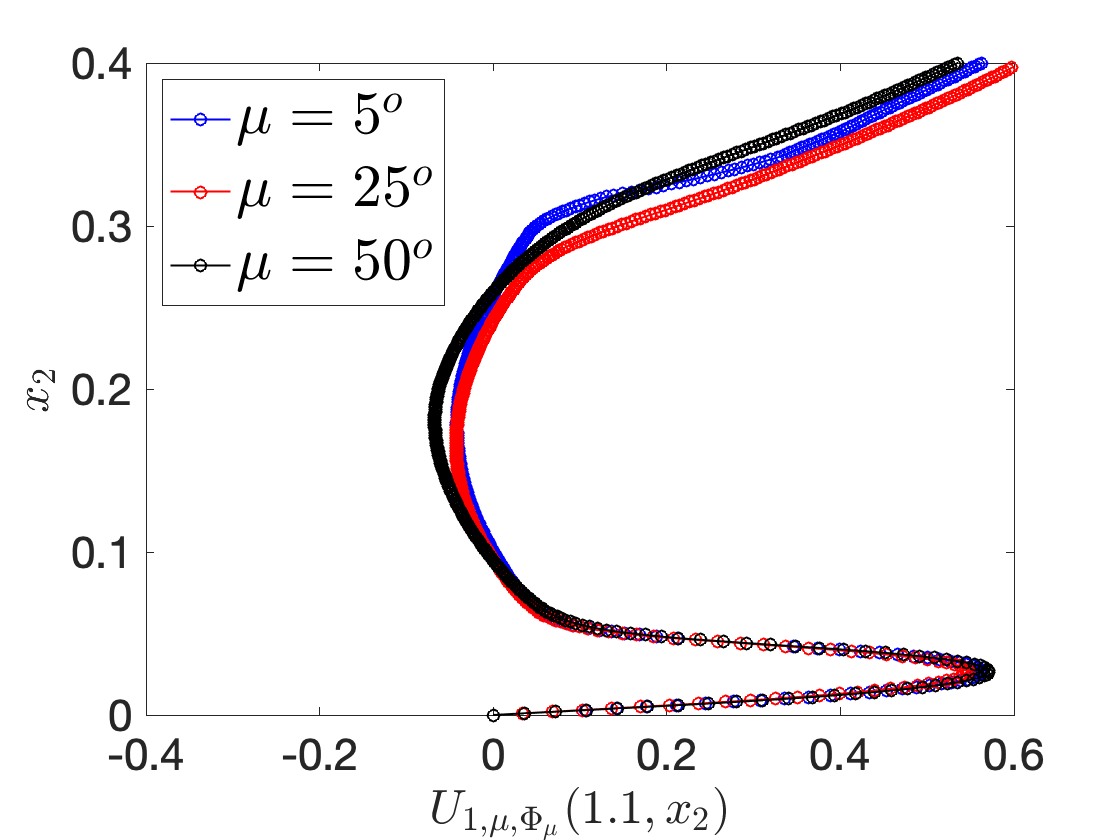}
  }
 \caption{Ahmed's body.
(a)  behavior of the horizontal velocity for $\mu=25^o$.
(b) behavior of the horizontal velocity profile for $x_1=1.1$ and various $x_2$ for three choices of $\mu$.
 }
\label{fig:ahmed_body_flow}
\end{figure}

Figure \ref{fig:ahmed_body_flow}(a) shows the behavior of $U_{1,\mu,\Phi_{\mu}}$ (horizontal velocity) for $\mu=25^o$, while Figure \ref{fig:ahmed_body_flow}(b) shows the profile of $U_{1,\Phi_{\mu}}$ for $x_1=1.1$ and $x_2\in [0,0.4]$ for three values of $\mu$: we observe the presence of  recirculating flow in the proximity of the body that becomes more prominent for  larger slant angles.

Since the boundary conditions are parameter-independent, we bypass the definition of the extension operator and we simply define ${\mathbf{e}}_{\mu}: ={\mathbf{U}}_{\bar{\mu}}^{\rm hf}$ with $\bar{\mu}=25^o$.
We equip the space ${\mathfrak{U}}$ with the inner product
\begin{equation}
\label{eq:inner_product_rans}
\langle {\mathbf{u}}, {\mathbf{v}} \rangle
\, = \,
\frac{1}{\lambda_{\rm v}} \,
\sum_{i=1,2} ( \mathbf{u}_i, \mathbf{v}_i )
+
\frac{1}{\lambda_{\rm p}} \,
\int_{\Omega} u_3 \, v_3 \, d{x}
+
\frac{1}{\lambda_{\rm \tilde{\nu}}} \,
( \mathbf{u}_4, \mathbf{v}_4 ),
\end{equation}
where $\lambda_{\rm v},\lambda_{\rm p},\lambda_{\rm \tilde{\nu}}>0$ are suitable positive constraints.
In this work, 
given the lifted snapshots
$\{ \mathring{{\mathbf{U}}}^k  \}_{k=1}^{n_{\rm train}}$, 
we define $\lambda_{\rm v},\lambda_{\rm p},\lambda_{\rm \tilde{\nu}}$ as the largest eigenvalues of the Gramian matrices
$\mathbf{C}^{\rm v},\mathbf{C}^{\rm p},\mathbf{C}^{\tilde{\nu}}$
such that
$$
\mathbf{C}_{i,j}^{\rm v}
=
\sum_{d=1,2} ( \mathring{\mathbf{U}}_d^{j}, 
\mathring{\mathbf{U}}_d^{i} ),
\quad
\mathbf{C}_{i,j}^{\rm p}
=
\int_{\Omega} \mathring{u}_3^{j} \, \mathring{u}_3^{i} \, d{x},
\quad
\mathbf{C}_{i,j}^{\tilde{\nu}}
=
( \mathring{\mathbf{u}}_4^{j}, 
\mathring{\mathbf{u}}_4^{i} ),
\quad
i,j=1,\ldots,n_{\rm train}.
$$
Our choice is motivated by the need to  properly ``weight'' variables 
(velocity, pressure, artificial viscosity) 
characterized by  different magnitudes (and different units!).
{
Due to our construction of the empirical test space, our choice of the inner product has also the effect of weighting the contributions of the  residuals associated with the different equations: similarly, 
in \cite[section IV]{washabaugh2016use} the authors consider a weighted residual to deal with the difference in magnitude  between the different  equations of the model.
}

\subsection{Numerical results}
\label{sec:RANS_numerics}
We train the reduced-order model based on $n_{\rm train}=40$ snapshots associated with equispaced parameters in $\mathcal{P}$. In order to generate the EQ weights $\boldsymbol{\rho}^{\rm eq}$, we consider additional $n_{\rm train,eq}=10$ randomly-selected parameters { (cf. comment after \eqref{eq:accuracy_constraint_nonlinear_b})}. The dual residual estimator is built using $tol_{\rm es}=10^{-4}$ and $tol_{\rm eq,r}=10^{-10}$. We assess performance based on $n_{\rm test}=10$ randomly-selected out-of-sample parameters.

Figure \ref{fig:performance_nonlinear}(a) shows the average out-of-sample error $E_{\rm avg}$ (cf. \eqref{eq:Eavg}) for various choices of the size $N$ of the reduced-order basis ${\mathbf{Z}}$ and for three different tolerances $tol_{\rm eq}$ for $\boldsymbol{\rho}^{\rm eq}$, the test space ${\mathbf{Y}}$ is chosen according to Algorithm \ref{alg:nonlinear_emp_test_space} with $J=2N$. As for the linear problem, we observe that the ROM guarantees near-optimal performance with respect to the projection error, for sufficiently tight tolerances. Figure \ref{fig:performance_nonlinear}(b) shows the { percentage} $Q/N_{\rm e} \cdot 100$ of sampled elements with respect to $N$ for various tolerances: for all test cases considered, $Q$ is less than $2\%$ of the size $N_{\rm e}$ of the complete mesh. { As for the linear case, speedups of the hyper-reduced ROM compared to the ROM with HF quadrature scale with $N_{\rm e}/Q$ and thus range from $\mathcal{O}(50)$ to $\mathcal{O}(300)$ for $N=1,\ldots,10$ and $tol_{\rm eq}=10^{-10}$.}
In Figure \ref{fig:residual_error_rans}, we replicate the results shown in Figure 
\ref{fig:residual_error} for the linear problem. For the test in Figure \ref{fig:residual_error_rans}(c), we consider various ROMs associated with $N=3,6,9$ and $tol_{\rm eq}=10^{-6},10^{-10},10^{-14}$. We observe that the dual residual is highly-correlated with the $H^1$ relative error: this empirical finding motivates the use of the (approximate) dual residual estimator to drive greedy sampling methods (cf. \cite{rozza2007reduced}).

\begin{figure}[H]
\centering
\subfloat[ ] {
\begin{tikzpicture}[scale=0.7]
\begin{semilogyaxis}[
xlabel = {\LARGE {$N$}},
  ylabel = {\LARGE {$E_{\rm avg}$}},
  line width=1.2pt,
  mark size=3.0pt,
  ymin=0.00001,   ymax=0.05,
  ]
  
\addplot[line width=1.pt,color=black,mark=o] table {data/RANS/Re3e3/proj.dat};
     \addplot[line width=1.pt,color=magenta,mark=square]  table {data/RANS/Re3e3/hf_quad.dat};
  \addplot[line width=1.pt,color=red,mark=diamond]  table {data/RANS/Re3e3/tolm6.dat}; 
  \addplot[line width=1.pt,color=blue,mark=triangle*] table {data/RANS/Re3e3/tolm10.dat}; 
    \addplot[line width=1.pt,color=violet,mark=otimes*] table {data/RANS/Re3e3/tolm14.dat}; 
    \label{ROMperformance:tol1em14};
    
\end{semilogyaxis}

\end{tikzpicture}
}
~~~
\subfloat[ ] {
\begin{tikzpicture}[scale=0.7]
\begin{axis}[
xlabel = {\LARGE {$N$}},
  ylabel = {\LARGE {$\%$ sampled elements}},
  line width=1.2pt,
  mark size=3.0pt,
 ymin=0,   ymax=2,
  ]
 
  \addplot[line width=1.pt,color=red,mark=diamond]  table {data/RANS/Re3e3/EQtolm6.dat};
  \addplot[line width=1.pt,color=blue,mark=triangle*] table {data/RANS/Re3e3/EQtolm10.dat};  
      \addplot[line width=1.pt,color=violet,mark=otimes*] table {data/RANS/Re3e3/EQtolm14.dat}; 
\end{axis}

\end{tikzpicture}
}

\bgroup
\sbox0{\ref{data}}%
\pgfmathparse{\ht0/1ex}%
\xdef\refsize{\pgfmathresult ex}%
\egroup
\caption[Caption in ToC]{
Ahmed's body problem.
Performance  with respect to $N$ and $tol_{\rm eq}$
($J=2 N$).
(a) average out-of-sample error $E_{\rm avg}$:
projection  \tikzref{ROMperformance:proj_error},
AMR ROM  with HF quadrature \tikzref{ROMperformance:hf_quad},
$tol_{\rm eq}=10^{-6}$ \tikzref{ROMperformance:tol1em6},
$tol_{\rm eq}=10^{-10}$ \tikzref{ROMperformance:tol1em10},
$tol_{\rm eq}=10^{-14}$ \tikzref{ROMperformance:tol1em14}
(b) percentage of sampled elements $Q/N_{\rm e}\cdot 100$ 
for three choices of $tol_{\rm eq}$ 
(see (a)).
}
 \label{fig:performance_nonlinear}
\end{figure}

\begin{figure}[h!]
\centering
\subfloat[ ] {
\begin{tikzpicture}[scale=0.7]
\begin{axis}[
xlabel = {\LARGE {$N$}},
  ylabel = {\LARGE {$J_{\rm r}$} },
  line width=1.2pt,
  mark size=3.0pt,
  ymin=1,   ymax=40
  ]
  \addplot[line width=1.pt,color=blue,mark=triangle*]  table {data/RANS/Re3e3/Jresm10.dat}; 
  \addplot[line width=1.pt,color=violet,mark=diamond] table {data/RANS/Re3e3/Jresm14.dat};
 
\end{axis}

\end{tikzpicture}
}
~~~
\subfloat[ ] {
\begin{tikzpicture}[scale=0.7]
\begin{axis}[
xlabel = {\LARGE  {$N$}},
  ylabel = {\LARGE {$\%$ sampled elements} },
  line width=1.2pt,
  mark size=3.0pt,
 ymin=0,   ymax=1,
  legend style={at={(0,0)},anchor=south west}
  ]
 
  \addplot[line width=1.pt,color=blue,mark=triangle*]  table {data/RANS/Re3e3/Qresm10.dat}; 
  \addplot[line width=1.pt,color=violet,mark=diamond] table {data/RANS/Re3e3/Qresm14.dat};
 
\end{axis}
\end{tikzpicture}
}

\subfloat[ ] {
\begin{tikzpicture}[scale=0.7]
\begin{loglogaxis}[
xlabel = {\LARGE {dual residual}},
  ylabel = {\LARGE {rel error}},
legend entries = {hf res vs error, estimated res vs error},
  line width=1.2pt,
  mark size=3.0pt,
 ymin=0.00001,   ymax=0.1,
legend style={at={(0.03,0.9)},anchor=west,font=\Large}
  ]
 
  \addplot[only marks,color=red,mark=square]  table {data/RANS/Re3e3/EQres.dat};
  
  \addplot[only marks,color=blue,mark=triangle*] table {data/RANS/Re3e3/hfres.dat};  
  
\end{loglogaxis}
\end{tikzpicture}
}
\caption{
Ahmed's body problem.
Dual residual norm estimation.
(a)  size of the empirical test space  ${J_{\rm r}}$ with respect to  $N$, for $tol_{\rm eq}=10^{-10}$ \tikzref{ROMperformance:tol1em10} and
$tol_{\rm eq}=10^{-14}$ \tikzref{ROMperformance:tol1em14}.
(b)  percentage of sampled elements $Q_{\rm r}/N_{\rm e}\cdot 100$     with respect to  $N$, for $tol_{\rm eq}=10^{-10}$ \tikzref{ROMperformance:tol1em10} and
$tol_{\rm eq}=10^{-14}$ \tikzref{ROMperformance:tol1em14}.
(c) exact  and estimated residual vs relative $H^1$ error.
($tol_{\rm eq,r}=10^{-10}, tol_{\rm es}=10^{-4}$, $J= 2N$)}
\label{fig:residual_error_rans}
\end{figure}  

\section{Conclusions}
\label{sec:conclusions}

In this work, we presented and numerically assessed a general approach for the treatment of parameterized geometries in projection-based pMOR. The approach relies on recently-developed EQ procedures to sample the HF mesh and ultimately reduce online costs. 
{ We discussed in detail the comparison with the more standard map-then-discretize approach;}
we presented the application to a linear diffusion problem to illustrate the key features of the methodology and analyze the offline and online computational and memory costs; then, we considered the extension to steady nonlinear PDEs and we discussed the application to the  two-dimensional RANS equations.

We here resorted to continuous FE discretizations;  { the extension to DG and finite volume discretizations is part of ongoing research.}
Furthermore,  we wish to combine the approach with adaptive sampling strategies to reduce offline costs: in this respect, the results in Figure \ref{fig:residual_error_rans} further motivate the use of residual-based error indicators to drive sampling. 
{ We also wish to investigate if the conservation properties of the hyper-reduced ROM derived in 
\cite{yano2019discontinuous,chan2020entropy} can be extended to problems with varying geometry,  within the  DtM framework.}
Finally, we shall combine the approach with solution-aware registration procedures, \cite{taddei2020registration,taddei2020space}, 
for the construction of the parameterized mapping ${\Phi}$.

\section*{Acknowledgements}
The authors acknowledge the support by European Union’s Horizon 2020 research and innovation programme under the Marie Skłodowska-Curie Actions, grant agreement 872442 (ARIA).
Tommaso Taddei also acknowledges the support of IdEx Bordeaux (projet EMERGENCE 2019).
The authors thank Dr. Andrea Ferrero (Politecnico di Torino) and Professor Angelo Iollo (Inria Bordeaux) for fruitful discussions on the numerical example of section \ref{sec:nonlinear_pb}.

 \appendix

\section{Explicit expressions of the RANS equations}
\label{sec:RANS_appendix}

We denote by ${u}$ the velocity field, by $p$ the fluid pressure, and by $\tilde{\nu}$ the turbulent viscosity; we further denote by 
${U} = [{u},p,\tilde{\nu}]$ the vector of state variables;  then, we introduce the RANS equations as
\begin{subequations}
\label{eq:RANS_strong}
\begin{equation}
\label{eq:RANS_stronga}
\left\{
\begin{array}{ll}
\displaystyle{
\partial_t {u} + {R}_{\rm mom} ( {U}  ) = {0} } &
{\rm in} \; \Omega\times \mathbb{R}_+ \\[3mm]
\nabla \cdot {u}   =0 
 & {\rm in} \; \Omega\times \mathbb{R}_+ \\[3mm]
 \displaystyle{
\partial_t \tilde{\nu} + {R}_{\nu} ( {U}  ) = 0 }& {\rm in} \; \Omega\times \mathbb{R}_+ \\ 
\end{array}
\right.
\end{equation}
where
\begin{equation}
\label{eq:RANS_strongb}
\left\{
\begin{array}{l}
 \displaystyle{
 {R}_{\rm mom} ( {U}  )
=
{u}\cdot \nabla {u} + \nabla p - \nabla \cdot {\tau},
\quad
{\tau}
= (\nu + \tilde{\nu}) \left(
\nabla {u} + \nabla {u}^T
\right)
}
\\[3mm]
 \displaystyle{
{R}_{\nu} ( {U}  )
=
{u}\cdot \nabla \tilde{\nu} -P + D
-\frac{3}{2} \left(
\nabla \cdot \left(  (\nu + \tilde{\nu}) \nabla \tilde{\nu}    \right)
+ 0.622 \| \nabla  \tilde{\nu}  \|_2^2
\right)
}
\\
\end{array}
\right.
\end{equation}
 where $\nu>0$ is the kinematic viscosity, $P, D$ are the wall production and destruction terms which are explicit functions of vorticity and turbulent viscosity (see \cite{Spalart1992, Allmaras2012}).  We are here interested in the limit solution for $t\to \infty$.
 \end{subequations}

FE discretizations of \eqref{eq:RANS_strong} can be obtained by projecting the equations onto finite dimensional FE spaces.  For stability reasons, we add SUPG and LSIC stabilization forms to the momentum equation
\begin{subequations}
\label{eq:RANS_FEM}
\begin{equation}
\left\{
\begin{array}{l}
 \displaystyle{
\mathcal{R}_{\rm mom}^{\rm supg}({U}, {v})
=
\sum_{k=1}^{N_{\rm e}} \int_{\texttt{D}_k} \;
\left(
\partial_t {u} +  {R}_{\rm mom} ( {U}  )
\right) \cdot 
\left(
\tau_{\rm supg}
 {u} \cdot \nabla {v}
\right) \, d {x}
}, \\[3mm]
 \displaystyle{
\mathcal{R}_{\rm mom}^{\rm lsic}({U}, {v})
=
\sum_{k=1}^{N_{\rm e}} \int_{\texttt{D}_k} \;
\tau_{\rm lsic}
\left(
\nabla\cdot  {u}
\right)  \,
\left( \nabla\cdot  {v} \right) \, d {x}
}, \\[3mm]
\end{array}
\right.
\end{equation}
where $\tau_{\rm supg}$ and $\tau_{\rm lsic}$ are defined as in 
\cite{Peterson2018}  and  \cite{braack2007stabilized}, respectively. We also add the PSPG term to the continuity equation
\begin{equation}
\mathcal{R}^{\rm pspg}({U}, q)
=
\sum_{k=1}^{N_{\rm e}} \int_{\texttt{D}_k} \;
\left(
\partial_t {u} +  {R}_{\rm mom} ( {U}  )
\right)  \,\cdot \, 
\left(   \tau_{\rm pspg}   \nabla q \right) \, d {x}
\end{equation}
with $\tau_{\rm pspg} =\tau_{\rm supg}$. Finally, we add  the SUPG  form  to the transport equation  for the turbulent viscosity:
\begin{equation}
\mathcal{R}_{\nu}^{\rm supg}({U}, \omega)
=
\sum_{k=1}^{N_{\rm e}} \int_{\texttt{D}_k} \;
\tau_{\rm supg}
\left(
\partial_t \tilde{\nu} +  {R}_{\nu} ( {U}  )
\right)  \,
\left(    {u} \cdot   \nabla \omega \right) \, d {x}
\end{equation}
\end{subequations}

In the implementation, we resort to pseudo-time integration to solve the HF system and we consider the strategy in \cite{Crivellini2013} to treat negative viscosities that might be generated during the iterations. Finally, boundary conditions are set as follows:
\begin{itemize}
\item
Inlet: ${u} = {u}_{\rm inlet}$,
$\tilde{\nu} = 3 \nu$.
\item
Top and outlet: homogeneous Neumann conditions.
\item
Bottom and body: ${u} = {0}$,
$\tilde{\nu} = 0$.
\end{itemize}
 
\section{Additional remarks}
\label{sec:appendix_remarks}

We provide below  further interpretations and remarks. 

\begin{remark}
\label{remark:dirichlet_norm}
\textbf{Choice of the boundary norm in \eqref{eq:bnd_norm}.}
Due to our choice of $\|  \cdot  \|_{\rm dir}$, the space $\mathfrak{H}$ is a discrete counterpart of $L^2(\Gamma_{\rm dir})$; for second-order elliptic PDEs, this choice is not consistent with the infinite-dimensional formulation, which involves the fractional Sobolev space $H^{1/2}(\Gamma_{\rm dir})$.
A consistent choice of $\|  \cdot  \|_{\rm dir}$ would be
$\|  \cdot  \|_{\rm dir} = \| \mathcal{H} \cdot   \|$: clearly, evaluations of the  latter require the solution to a HF problem of size $N_{\rm hf}$.  Our choice is thus a compromise between mathematical rigor and simplicity of implementation.
\end{remark}

\begin{remark}
\label{remark:RID}
\textbf{Reduced integration domain.}
Following \cite{ryckelynck2009hyper,fritzen2018algorithmic}, we might refer to the region $ \bigcup_{k \in \texttt{I}_{\rm eq}}   \texttt{D}_k $ as to \emph{reduced integration domain} (RID). 
Note that, while 
the EQ procedure employed in this paper directly reduces the integration domain without modifying the test functions, 
the hyper-reduction approach in 
\cite{fritzen2018algorithmic} 
--- 
similarly to \cite{chaturantabut2010nonlinear}
---
ultimately performs projection of the residual on the nodes of the mesh: as a result, 
for the linear problem\footnote{The same result holds for the nonlinear case.} 
considered in section \ref{sec:linear_case},
 the latter enforces that
 $$
 {\rm if} \; 
 \mathcal{R}_{\mu} \left( \mathbf{Z} \boldsymbol{\alpha} + \widehat{\mathbf{e}}_{\mu}, {e}_i, {\Phi}_{\mu}  \right) = 0
 \;\;
 \forall \, i=1,\ldots,N_{\rm hf}
 \;\;
 {\rm then \;}
 \boldsymbol{\alpha} = \widehat{\boldsymbol{\alpha}}_{\mu}, 
  $$
which is in general false for the EQ approach considered here.
We further remark that the strategy employed in \cite{ryckelynck2009hyper,fritzen2018algorithmic}  to construct the RID is different, and the authors do not explicitly discuss the problem of parameterized geometries.
\end{remark}

\begin{remark}
\textbf{Connection between linear and nonlinear procedures.}
From the discussion in section \ref{sec:hf_model}, we deduce that 
the FE solver for \eqref{eq:variational_formulation}  shall include a routine of the form:
$$
\left[
\mathbf{F}_{k,\Phi}^{\rm un},  \mathbf{A}_{k,\Phi}^{\rm un}
\right]
=
\texttt{local}{\_}\texttt{assembler} \, \left(
{\Phi} ( \texttt{X}_k^{\rm hf}  ), \; \texttt{ref}, \;  \texttt{input}{\_}\texttt{data} \right),
$$
where $\mathbf{F}_{k,\Phi}^{\rm un},  \mathbf{A}_{k,\Phi}^{\rm un}$ are the contributions to the global system associated with the $k$-th element,
${\Phi} ( \texttt{X}_k^{\rm hf}  )$ are the nodes of the deformed mesh associated with the $k$-th element, 
$\texttt{ref}$ and $\texttt{input}{\_}\texttt{data}$ are data structures introduced in section \ref{sec:hf_solver_nonlinear} --- note that $\mathbf{F}_{k,\Phi}^{\rm un} \equiv \mathbf{0}$ in the example of section \ref{sec:linear_case}.
Then, recalling 
\eqref{eq:unassembled_reduced_basis}, we find that 
\eqref{eq:GalROMb} could be rewritten as
$$
\left\{
\begin{array}{ll}
\displaystyle{
\left( \widehat{\mathbf{A}}_{\mu}^{\rm eq}  \right)_{n,n'} =
\sum_{k \in \texttt{I}_{\rm eq}  }  \; 
\rho_k^{\rm eq} \; 
\sum_{i,i'=1}^{n_{\rm lp}} 
\mathbf{Z}_{i,k,n}^{\rm un}
\;
\left(  \mathbf{A}_{k,\Phi}^{\rm un} \right)_{i,i'} 
\mathbf{Z}_{i',k,n'}^{\rm un},
} &\\[3mm]
\displaystyle{
\left( \widehat{\mathbf{F}}_{\mu}^{\rm eq}  \right)_n=
\sum_{k \in \texttt{I}_{\rm eq}  }  \; 
\rho_k^{\rm eq} \; 
\left(
\sum_{i=1}^{n_{\rm lp}} 
\mathbf{Z}_{i,k,n}^{\rm un}
\left(
\left( \mathbf{F}_{k,\Phi}^{\rm un} \right)_{i}
-
\sum_{i'=1}^{n_{\rm lp}} 
\left(  \mathbf{A}_{k,\Phi}^{\rm un} \right)_{i,i'} 
\left( \widehat{\mathbf{e}}^{\rm un}  \right)_{i',k} 
\right)
\right)
} &
\\
\end{array}
\right.
$$
The latter identity 
is the linear (and scalar) counterpart of the general relation provided in \eqref{eq:assembler_ROM_nonlinear}.
\end{remark}

\begin{remark}
\label{remark:positivity}
\textbf{Positivity of the empirical weights.}
Positivity of the weights $\boldsymbol{\rho}^{\rm eq}$ in \eqref{eq:weighted_residual} is justified by the analysis in \cite{yano2019discontinuous}; on the other  hand, the constraint $\boldsymbol{\rho}^{\rm eq,r} \geq \mathbf{0}$ is not strictly necessary and might be omitted.
We might thus construct the weights 
$\boldsymbol{\rho}^{\rm eq,r}$ by approximating 
the  solution to the sparse representation problem 
\begin{equation}
\label{eq:no_positivity_constraint}
\min_{  \boldsymbol{\rho} \in \mathbb{R}^{N_{\rm e}} }
\;
\| \boldsymbol{\rho}   \|_{\ell^0},
\quad
{\rm s.t.}  \; \; \|\mathbf{G}^{\rm r} \boldsymbol{\rho} - \mathbf{b}^{\rm r}   \|_{\star} \leq \delta.
\end{equation}
To apply nonnegative least-squares  to \eqref{eq:no_positivity_constraint}, we first observe that (see the discussion in \cite[section 2.3.1]{taddei2019offline}), given  any solution 
$(\boldsymbol{\rho}_+^{\rm eq,r}, \boldsymbol{\rho}_-^{\rm eq,r} )$ to 
\begin{equation}
\label{eq:no_positivity_constraint_alter}
\min_{  \boldsymbol{\rho}_+,\boldsymbol{\rho}_-  \in \mathbb{R}^{N_{\rm e}} }
\;
\| \boldsymbol{\rho}_+   \|_{\ell^0}+\| \boldsymbol{\rho}_-   \|_{\ell^0},
\quad
{\rm s.t.}  \; \; \|\mathbf{G}^{\rm r} (\boldsymbol{\rho}_+ - \boldsymbol{\rho}_- )  - \mathbf{b}^{\rm r}   \|_{\star} \leq \delta,
\;\;
\boldsymbol{\rho}_+,\boldsymbol{\rho}_- \geq \mathbf{0},
\end{equation}
 $\boldsymbol{\rho}^{\rm eq,r} =\boldsymbol{\rho}_+^{\rm eq,r}-\boldsymbol{\rho}_-^{\rm eq,r} $  solves \eqref{eq:no_positivity_constraint}. Problem \eqref{eq:no_positivity_constraint_alter} can be tackled using nonnegative least-squares; note, however, that \eqref{eq:no_positivity_constraint_alter} has twice  as many degrees of freedom as \eqref{eq:sparse_representation_res} and is thus more challenging to be approximated. 
\end{remark}

\begin{remark}
\label{remark:justification_residual}
\textbf{Mathematical justification of the choice of $\mathfrak{Y}_{J_{\rm r}}^{\rm r}, \boldsymbol{\rho}^{\rm eq,r}$ in section \ref{sec:a_posteriori}.}
We denote by $\Pi_{\mathfrak{Y}_{J_{\rm r}}^{\rm r}}: \mathfrak{U} \to \mathfrak{Y}_{J_{\rm r}}^{\rm r}$ the orthogonal projection operator onto $\mathfrak{Y}_{J_{\rm r}}^{\rm r}$. 
We can then exploit the same argument as in 
\cite[Proposition 2]{taddei2019offline} to obtain
\begin{equation}
\label{eq:error_estimate_residual}
\big|
\widehat{\mathfrak{R}}_{\mu}^{\rm eq}
-
{\mathfrak{R}}_{\mu}^{\rm hf}(  \widehat{\mathbf{u}}_{\mu})
\big|
\; \leq \; 
\frac{1}{
{\mathfrak{R}}_{\mu}^{\rm hf}(  \widehat{\mathbf{u}}_{\mu}) +
\| \widehat{\mathbf{R}}_{\mu}^{\rm hf, r}  \|_2 }
\underbrace{
\| \boldsymbol{\psi}_{\mu}^{\rm r} -  
\Pi_{\mathfrak{Y}_{J_{\rm r}}^{\rm r}} \boldsymbol{\psi}_{\mu}^{\rm r}   \|^2 }_{\rm (I)}
+
\underbrace{
\| \widehat{\mathbf{R}}_{\mu}^{\rm eq,r} - \widehat{\mathbf{R}}_{\mu}^{\rm hf,r}   \|_2 }_{\rm (II)}.
\end{equation}
We observe that (I) is associated with the approximation properties of the empirical test space for the manifold $\mathfrak{M}_{\rm test}^{\rm r}$, while (II) --- which is \eqref{eq:accuracy_constraint_res} --- is related to the accuracy of the quadrature rule.
\end{remark}

\bibliographystyle{abbrv}
\bibliography{all_refs}
 
\end{document}